\documentclass{amsart}

\title{On the $K$-theory of groups with
                   finite asymptotic dimension}
\author{Arthur Bartels and David Rosenthal}

\date{June 12, 2007}

\thanks{This work was supported by the SFB 478
 ``Geometrische Strukturen in der Mathematik''.}

\usepackage{hyperref}
\usepackage{enumerate,amssymb}
\usepackage[arrow,curve,matrix,tips]{xy}
  \SelectTips{eu}{10} \UseTips
\usepackage{amscd}
\usepackage[mathscr]{eucal}
\usepackage{amssymb}
\usepackage{calc}
\usepackage{bibunits}

\DeclareMathAlphabet{\matheurm}{U}{eur}{m}{n}

\newcommand{\pt}{\text{pt}}

\newcommand{\Or}{\matheurm{Or}}

\DeclareMathOperator{\id}{id}

\DeclareMathOperator{\supp}{supp}

\newcommand{\Fin}{{\mathcal{F}in}}

\newcommand{\Ebar}{\underbar{\rm{E}}\Gamma}

  \newcommand{\IK}{\mathbb{K}}
  \newcommand{\IL}{\mathbb{L}}
  
  \newcommand{\IN}{\mathbb{N}}

  \newcommand{\IR}{\mathbb{R}}

  \newcommand{\IZ}{\mathbb{Z}}

  \newcommand{\cala}{\mathcal{A}}
  \newcommand{\calb}{\mathcal{B}}
  \newcommand{\calc}{\mathcal{C}}

  \newcommand{\calf}{\mathcal{F}}

  \newcommand{\cals}{\mathcal{S}}
  
  \newcommand{\calu}{\mathcal{U}}

  \newcommand{\calx}{\mathcal{X}}
  \newcommand{\caly}{\mathcal{Y}}

  \newcommand{\bfA}{{\mathbf A}}
  \newcommand{\bfB}{{\mathbf B}}

  \newcommand{\bfK}{{\mathbf K}}
  \newcommand{\bfL}{{\mathbf L}}

  \newcommand{\bfR}{{\mathbf R}}
  \newcommand{\bfS}{{\mathbf S}}
  \newcommand{\bfT}{{\mathbf T}}

\newcommand{\x}{{\times}}

\newcommand{\dd}{{\partial}}
\newcommand{\pl}{{\mathit{pl}}}
\newcommand{\onenorm}[1]{{\left\| #1 \right\|_1}}
\newcommand{\twonorm}[1]{{\left\| #1 \right\|_2}}

\newcommand{\lieg}{{\mathbf g}}
\newcommand{\lieh}{{\mathbf h}}

\newcommand{\liel}{{\mathbf l}}
\newcommand{\liev}{{\mathbf v}}

\DeclareMathOperator{\Ad}{Ad}

\newenvironment{numberlist}
  {\begin{list}{}%
   {%
    \setlength{\leftmargin}{\labelwidth+\labelsep}%
   }%
  }%
  {\end{list}}

\DeclareMathOperator{\Map}{Map}

\theoremstyle{plain}
\newtheorem{theorem}{Theorem}[section]

\newtheorem{lemma}[theorem]{Lemma}

\newtheorem{proposition}[theorem]{Proposition}

\newtheorem*{theorema}{Theorem A}
\newtheorem*{theoremb}{Theorem B}
\newtheorem*{corollary*}{Corollary}

\theoremstyle{definition}
\newtheorem{definition}[theorem]{Definition}

\theoremstyle{remark}
\newtheorem{remark}[theorem]{Remark}
\newtheorem{notation}[theorem]{Notation}

\makeatletter\let\c@equation=\c@theorem\makeatother

\begin{document}

\newlength{\origlabelwidth}
\setlength\origlabelwidth\labelwidth

\begin{abstract}
It is proved that the assembly maps in algebraic $K$-
and $L$-theory with respect to the family of finite
subgroups is injective
for groups $\Gamma$ with finite asymptotic dimension that admit
a finite model for $\Ebar$.
The result also applies to certain groups that admit only
a finite dimensional model for $\Ebar$.
In particular, it applies to discrete subgroups of
virtually connected Lie groups.
\end{abstract}

\maketitle


\begin{bibunit}

\section*{Introduction}

Assembly maps in algebraic $K$- and $L$-theory
are designed to
study the $K$- and $L$-theory of group rings $R[\Gamma]$,
which contain important geometric information about
manifolds with fundamental group $\Gamma$.
Similarly, the Baum-Connes map is used to analyze the
topological $K$-theory of the reduced $C^*$-algebra of $\Gamma$.
An important class of groups that has been
studied in recent years,
by topologists and analysts alike,
is the class of discrete groups
with {\it finite asymptotic dimension}.
Yu proved the Novikov conjecture for groups with
finite asymptotic
dimension that admit a finite classifying
space~\cite{Yu-finite-asymptotic-dimension-Novikov}.
He achieved this by using controlled techniques to prove the
injectivity of the Baum-Connes map for such groups.
Later on, Higson was able to remove the finite classifying
space assumption~\cite{Higson-Bivariant-K-theory}.
In~\cite{Bartels-finite-asymp-dim}, the first author proved
the algebraic $K$- and $L$-theory versions of Yu's work
by developing a squeezing theorem for higher
algebraic $K$-theory
to use with the approach established
in~\cite{Yu-finite-asymptotic-dimension-Novikov}.
For algebraic $K$-theory, this was also achieved in
\cite{Carlsson-Goldfarb-finite-asymptotic}.

The purpose of this paper is to extend the results
from~\cite{Bartels-finite-asymp-dim} by relaxing the
finite ${\rm B}\Gamma$ assumption to allow for groups $\Gamma$
with torsion. Specifically, we prove the following theorem.

\begin{theorema}
\label{thm:main-K-theory}
Let $\Gamma$ be a discrete group and $R$ a ring.
Assume that there is a finite dimensional $\Gamma$-$CW$-model
for the universal space for proper $\Gamma$-actions,
$\Ebar$, and
assume that there is a $\Gamma$-invariant metric
on $\Ebar$ such that
$\Ebar$ is uniformly $\Fin$-contractible,
is a complete proper path metric space
and has finite asymptotic dimension.
Then the assembly map,
\begin{equation}
\label{K-assembly-map-ring}
H_*^{\Gamma} (\Ebar;
  \bfK_R) \to K_*(R[\Gamma]),
\end{equation}
in algebraic $K$-theory, is a split injection.
\end{theorema}

The notion of  uniform $\Fin$-contractibility is a
strengthening of uniform contractibility and is
introduced in Definition~\ref{def:fin-contractible}.
The asymptotic dimension of a metric space was introduced
by Gromov \cite{Gromov-asymptotic-invariants}
and is reviewed in Section~\ref{sec:asymptotic-dimension}.

Note that Theorem~A\ref{thm:main-K-theory}
has interesting consequences for Whitehead groups.
The classical assembly map
\begin{equation}
\label{K-assembly-map-classical}
H_n(B\Gamma;\bfK_R) \to K_n(R[\Gamma])
\end{equation}
considered in
\cite{Bartels-finite-asymp-dim,
           Carlsson-Goldfarb-finite-asymptotic},
factors through the assembly map
\eqref{K-assembly-map-ring} via
\begin{equation}
\label{relative-assembly}
H_n(B\Gamma;\bfK_R) \cong
             H_n^\Gamma(E\Gamma;\bfK_R) \to
             H_n^\Gamma(\Ebar;\bfK_R).
\end{equation}
In particular, the cokernel of \eqref{K-assembly-map-classical}
contains the cokernel of \eqref{relative-assembly}
in the situation of
Theorem~A\ref{thm:main-K-theory}.
This cokernel can be evaluated using an Atiyah-Hirzebruch
spectral sequence
(see \cite[Theorem~4.7]{Davis-Lueck-assembly}),
or computed rationally using the equivariant
Chern character from
\cite[Theorem~0.3]{Lueck-Chern-character},
\cite[Theorem~173]{Lueck-Reich-survey}.
If $R = \IZ$ and $n= 1$, then the cokernel
of \eqref{K-assembly-map-classical} is the Whitehead group.
In this way, Theorem~A\ref{thm:main-K-theory}
implies non-vanishing results for the Whitehead group.

In Section~\ref{sec:subgroups-of-liegroups}
it is shown that the assumptions of
Theorem~A\ref{thm:main-K-theory}
are satisfied for all discrete subgroups of Lie groups with a
finite number of components.
Also, if there is a cocompact $\Gamma$-$CW$-model for $\Ebar$,
then, with any $\Gamma$-invariant metric, it will be
quasi-isometric to $\Gamma$ when $\Gamma$ is equipped
with a word length metric.
By Lemma~\ref{lem:finite-EG-is-uniformly-Fin-contractible},
such a model is uniformly $\Fin$-contractible.
Thus, Theorem~A\ref{thm:main-K-theory}
implies the following corollary.

\begin{corollary*}
Let $\Gamma$ be a group and assume that one of the
following two conditions is satisfied:
\begin{enumerate}
\item $\Gamma$ is a discrete subgroup of a virtually connected
      Lie group.
\item $\Gamma$ has finite asymptotic dimension and admits
      a cocompact $\Gamma$-$CW$-model for $\Ebar$.
\end{enumerate}
Then the assembly map \eqref{K-assembly-map-ring}
in algebraic $K$-theory
is split injective for every ring $R$.
\end{corollary*}

The techniques used to prove Theorem~A\ref{thm:main-K-theory}
allow for a very similar proof of the corresponding result in
$L$-theory if {\em ultimate lower quadratic}
$L$-theory, $L_*^{\langle -\infty \rangle}$, is used.
The only difference is that
the compatibility of $L$-theory with infinite products
is only known under an
additional $K$-theory assumption.
This forces the extra hypothesis in the following theorem.

\begin{theoremb}
\label{thm:main-L-theory}
Let $\Gamma$ be a discrete group and $R$ a ring with involution.
Assume that there is a finite dimensional $\Gamma$-$CW$-model
for the universal space for proper $\Gamma$-actions, $\Ebar$,
and assume that there is a $\Gamma$-invariant
metric on $\Ebar$ such
that $\Ebar$ is uniformly $\Fin$-contractible,
is a complete proper path metric space
and has finite asymptotic dimension.
Further assume that for every finite subgroup $G$ of $\Gamma$,
the group $K_{-i}(R[G])$ vanishes for sufficiently large $i$.
Then the assembly map,
\begin{equation}
\label{L-assembly-map-ring}
H_*^{\Gamma} (\Ebar;
  \bfL_R) \to L^{\langle -\infty \rangle}_*(R[\Gamma]),
\end{equation}
in $L$-theory, is a split injection.
\end{theoremb}

As in the $K$-theory case,
Theorem~B\ref{thm:main-L-theory} can be used to obtain
non-vanishing results for the cokernel of the
assembly map
\[
H_*(B\Gamma;\bfL_R) \to L_*^{\langle -\infty \rangle}(R[\Gamma]).
\]
It is well-known that the Novikov conjecture
on the homotopy invariance of higher signatures
is implied by the rational injectivity of this map.
Clearly, Theorem~B\ref{thm:main-L-theory}
also implies the following result.

\begin{corollary*}
Let $\Gamma$ be a group and assume that one of the
following two conditions is satisfied:
\begin{enumerate}
\item $\Gamma$ is a discrete subgroup of a virtually connected
      Lie group.
\item $\Gamma$ has finite asymptotic dimension and admits
      a cocompact $\Gamma$-$CW$-model for $\Ebar$.
\end{enumerate}
Then the assembly map \eqref{L-assembly-map-ring} in $L$-theory
is split injective for every ring $R$ with involution such
that for every finite subgroup $G$
of $\Gamma$, the group $K_{-i}(R[G])$ vanishes for
sufficiently large $i$.
\end{corollary*}

It is interesting to compare the finiteness conditions in
Theorems~A\ref{thm:main-K-theory} and B\ref{thm:main-L-theory}
with the finiteness assumption in the rational
injectivity result
of B\"okstedt-Hsiang-Madsen
\cite{Boekstedt-Hsiang-Madsen-cyclotomic-trace}
for the assembly map~\eqref{K-assembly-map-classical} (with
$R = \IZ$),
where it is assumed that the integral homology,
$H_n(B\Gamma;\IZ)$,
is finitely generated for every $n$.
The only other injectivity results that we are aware of
that apply to
all subgroups of virtually connected Lie groups are
those for the Baum-Connes assembly map
by Kasparov~\cite{Kasparov-equivariant-KK-and-Novikov}
and Higson~\cite{Higson-Bivariant-K-theory}
(compare \cite[Section~4]{Higson-Roe-amenable-actions+Novikov}).
This also implies the rational injectivity of
\eqref{L-assembly-map-ring} for $R = \IZ$.
It should also be noted that
Ferry-Weinberger
\cite{Ferry-Weinberger-Cuvature-tangentiality+controlled}
prove the Novikov conjecture and
the rational injectivity of~\eqref{K-assembly-map-classical}
(with $R = \IZ$) for
fundamental groups of non-positively curved manifolds
that are not necessarily compact.
By \cite{Bartels-on-the-domain}, our results are in accordance
with the Farrell-Jones conjecture
\cite{Farrell-Jones-Isomorphism-Conjectures}.
For more information about the Baum-Connes and Farrell-Jones
conjectures we recommend \cite{Lueck-Reich-survey}.

In Sections~\ref{sec:proof} and~\ref{sec:L-theory},
Theorems~\ref{thm:split-injectivity-with-coeficients}
and \ref{thm:split-injectivity-with-coeficients-L-theory}
are proven which are
slightly stronger than the theorems stated in this introduction.
There we will be considering assembly maps for
an additive category
with a $\Gamma$-action, as introduced in
\cite{BR-coefficients}.
For example, this more general setup enables
one to study the $K$- and $L$-theory of crossed product
rings (see \cite[Section~6]{BR-coefficients}).
Theorems~A\ref{thm:main-K-theory} and
B\ref{thm:main-L-theory} follow from
Theorems~\ref{thm:split-injectivity-with-coeficients}
and \ref{thm:split-injectivity-with-coeficients-L-theory}
respectively,
by using the category of finitely generated free
$R$-modules with the trivial action of $\Gamma$.

The authors would like to thank Burkhard Wilking
for his help with
Section~\ref{sec:subgroups-of-liegroups}.
We thank Lizhen Ji for pointing out an incorrect statement
in a previous version of this paper.


\section{Uniform contractibility}

Let $A$ be a subset of a metric space $X$, and let $R > 0$.
Denote by $A^R$ the set of all points $x$ in $X$
for which $d(x,A) \leq R$.
If $A = \{ x \}$ consists of just one point, then we will
abbreviate $x^R = \{ x \}^R$.

\begin{definition}
[Uniformly $\Fin$-contractible]
\label{def:fin-contractible}
Let $\Gamma$ act by isometries on a metric space $X$.
Then $X$ is said to be {\it uniformly $\Fin$-contractible} if
for every finite subgroup $G$ of $\Gamma$ and every $R > 0$
there is an $S > 0$ such that the following holds:
If $B$ is a $G$-invariant subset of $X$ of diameter less
than $R$, then the
inclusion, $B \to B^S$, is $G$-equivariantly null homotopic.
In particular, $B^S$ contains a fixed point for the
action of $G$.
\end{definition}

\begin{remark}
[Uniformly contractible]
If $\Gamma$ is the trivial group, or more generally if $\Gamma$
is torsion-free,
then $X$ is uniformly $\Fin$-contractible if and only if
$X$ is {\em uniformly contractible}. That is, for
every $R > 0$ there is an $S > R$ such that for every $x \in X$,
the inclusion, $x^R \to x^S$, is null homotopic.
\end{remark}

Let $X$ be a metric space with an isometric action
of a finite group $G$. Let $q \colon X \to X/G$ denote
the quotient map. Then
\[
d_{X/G} (y,y') = d_X ( q^{-1}(y), q^{-1}(y') ),
                                   \quad y,y' \in X/G
\]
defines a metric on $X/G$.
We will always consider quotients of metric spaces by
an isometric
action of a finite group as metric spaces using this metric.

\begin{notation}
\label{nota:X^cals}
Let $X$ be a space with an action of a group $G$, and let
$\cals$ be a collection of subgroups of $G$.
Define $X^\cals$ to be the union of all fixed sets
$X^H$, where $H$ is in $\cals$.
If $\cals=\{ G \}$, then $X^\cals$ is simply $X^G$.
If $\cals$ is closed under conjugation by elements of $G$,
then $X^\cals$ is a $G$-invariant subspace of $X$.
\end{notation}

\begin{lemma}
\label{lem:uniform-contractible-quotients}
Let $X$ be a metric space that is uniformly $\Fin$-contractible
with respect to an isometric action of a group $\Gamma$.
If $G$ is a finite subgroup of $\Gamma$ and $\cals$ is a
collection of subgroups of $G$ that is closed under
conjugation by $G$, then the quotient
$X^\cals / G$ is uniformly contractible.
\end{lemma}

\begin{proof}
Let $R > 0$ be given.
Since $X$ is assumed to be uniformly $\Fin$-contractible
and $G$ is finite, there
is an $S > 0$ such that for every subgroup $H$ of $G$ and
every $H$-invariant subset $B \subset X$
of diameter less than or equal to $2R|G|$,
the inclusion, $B \to B^S$, is $H$-equivariantly null homotopic.

Let $y \in X^\cals / G$ and $x \in q^{-1}(y)$, where
$q \colon X^\cals \to X^\cals / G$ is the quotient map.
Let $H$ be the subgroup of $G$ consisting of all $g \in G$
for which there are $g_1, \dots, g_n \in G$ such that
$g_1 = e$, $g_n = g$ and $d(g_i x,g_{i+1} x) \leq 2R$, for
$i=1,\dots,n-1$.
Then the diameter of $B = Hx^R \subset X$ is
bounded by $2R|G|$.
Therefore, the inclusion $B \to B^S$ is
$H$-equivariantly null homotopic.
In particular, there is a point
$z \in B^S$ that is fixed by $H$.
For $g \in G - H$, $gB \cap B = \emptyset$.
Therefore, the inclusion, $G\cdot B \to G\cdot B^S$, is
$G$-equivariantly
homotopic to a map that sends $gB$ to $gz$.
By $G$-equivariance, this homotopy can be restricted
to $X^\cals$, which induces the required null homotopy
on the quotient.
\end{proof}

\begin{lemma}
\label{lem:finite-EG-is-uniformly-Fin-contractible}
Let $\Gamma$ be a group such that there is a finite
$\Gamma$-$CW$-model, $X$, for $\Ebar$.
Let $d$ be a $\Gamma$-invariant metric on $X$.
Then $X$ is uniformly $\Fin$-contractible.
\end{lemma}

\begin{proof}
Let $X$ be a finite $\Gamma$-$CW$-model for $\Ebar$.
That is, there is a proper $\Gamma$-$CW$-complex $X$ with
finitely many $\Gamma$-cells
such that $X^G$ is
contractible for every finite subgroup $G$ of
$\Gamma$ and empty otherwise.
Note that $X$ is a locally compact space since
the $\Gamma$-action is cocompact and proper, and
for every finite subgroup $G$ of $\Gamma$, $X$ is
$G$-equivariantly contractible.

Let $R > 0$ be given.
If $B$ is a finite $G$-invariant subcomplex of
$X$ of diameter less than $R$, then there is an $S = S(B, G)$
such that $B \to B^S$ is $G$-equivariantly null homotopic.
If $\gamma \in \Gamma$, then
$\gamma B \to \gamma B^S$ is
$G^\gamma = \gamma G \gamma^{-1}$-equivariantly null homotopic
since the metric is $\Gamma$-invariant.
Consider the set of all pairs $(B,G)$, where $G$ is a finite
subgroup of $\Gamma$ and $B$ is a $G$-invariant
subcomplex of $X$
whose diameter is less than $R$.
On this set, $\gamma \in \Gamma$ acts by sending
$(B,G) \mapsto (\gamma B, G^\gamma)$.
Since the $\Gamma$-action on $X$ is proper,
the quotient by this action is finite.
Therefore, we can choose $S$ independent of $B$.
(In fact, $S$ can be chosen independent of both $B$ and $G$.)
\end{proof}


\section{Asymptotic dimension}
\label{sec:asymptotic-dimension}

Let $\calu$ be an open cover of the metric space $X$.
The cover $\calu$ is called {\em locally finite}
if every compact subset of $X$ meets only finitely many
members of $\calu$.
The dimension of an open cover $\calu$ is defined to be the
smallest number $n$ such that each $x$ in $X$ is contained in
at most $n+1$ members of $\calu$.
This is also the dimension
of the associated simplicial complex $|\calu|$.
If the diameters of the open sets in $\calu$ are
uniformly bounded,
then $\calu$ will be called a {\em bounded cover}.
The {\it asymptotic dimension}
\cite[p.28]{Gromov-asymptotic-invariants}
of $X$ is the smallest integer $n$ such that for
any $R > 0$, there exists an $n$-dimensional bounded
cover $\calu$ of $X$ whose Lebesgue number is at least $R$.

\begin{lemma}
\label{lem:finite-asymp-dim-and-locally-finite}
Let $X$ be a proper metric space of finite
asymptotic dimension $n$.
Then for every $\alpha$ there exists a locally finite
$n$-dimensional bounded
cover $\calu$ of $X$ whose Lebesgue number is at least $\alpha$.
\end{lemma}

\begin{proof}
Let $\calu$ be an $n$-dimensional bounded cover
of $X$ such that for every $x \in X$ there is $U \in \calu$
that contains the closed ball $x^R$.
For each $U \in \calu$, let $U^{-\alpha}$ be the open set
consisting of all points $x \in X$ for
which $x^\alpha \subset U$.
Then $\calu^{-\alpha} = \{ U^{-\alpha} \; | \; U \in \calu \}$
is an open cover of $X$.
Since $X$ is a proper metric space, we can find a
subcollection $\calu_0 \subseteq \calu$ and an open set
$U' \subset U^{-\alpha}$ for
each $U \in \calu_0$
such that $\{ U' \; | \; U \in \calu_0 \}$
is a locally finite set.
For each $U \in \calu_0$, let $U''$ be the interior
of $(U')^\alpha$.
The collection $\calu'' = \{ U'' \; | \; U \in \calu \}$
is also locally finite.
Since $U'' \subset U$ for $U \in \calu_0$,
$\calu''$ is a bounded cover of $X$ of dimension at most $n$.
By construction, the Lebesgue number of $\calu''$ is at least
$\alpha$.
\end{proof}

Recall the definition of $X^\cals$ in Notation~\ref{nota:X^cals}.

\begin{lemma}
\label{lem:asymptotic-dimension-and-finite-quotients}
Let $X$ be a metric space with a proper isometric action
of the group $\Gamma$.
Let $G$ be a finite subgroup of $\Gamma$.
Let $\cals$ be a  collection of  subgroups
of $G$ that is closed under conjugation by $G$.
If $X$ has finite asymptotic dimension, then
the quotient $X^\cals / G$ has finite asymptotic
dimension.
\end{lemma}

\begin{proof}
Let $n$ be the asymptotic dimension of $X$, and let $R > 0$
be given.
Because the asymptotic dimension of a subspace is bounded by
the asymptotic dimension of the ambient space, there exists
an $n$-dimensional bounded cover $\calu$ of $X^\cals$ whose
Lebesgue number is at least $R$.
Let $p \colon X^\cals \to X^\cals / G$ be the quotient map.
Define $p_* \calu = \{ p(U) \; | \; U \in \calu \}$.
It is easy to check that $p_* \calu$ is a bounded cover
whose Lebesgue number is at least $R$.
For every $x \in X^\cals / G$, $p^{-1}(x)$ contains no more
than $|G|$ points. Therefore, the dimension of $p_* \calu$
is at most $(n+1)|G| - 1$.
\end{proof}


\section{Discrete subgroups of Lie groups}
\label{sec:subgroups-of-liegroups}

In this section it is proven that discrete subgroups of
virtually connected Lie groups satisfy the assumptions
of Theorems~\ref{thm:split-injectivity-with-coeficients}
and \ref{thm:split-injectivity-with-coeficients-L-theory}.
Let $\Gamma$ be a discrete subgroup of a virtually
connected Lie group $G$.
If $K$ is a maximal compact subgroup of $G$, then $G/K$
is a finite dimensional
$\Gamma$-CW complex that is a model for the universal
proper $\Gamma$-space $\Ebar$
\cite[Theorem~4.4]{Lueck-survey-classifying-spaces}.
Furthermore, there exists a $G$-invariant Riemannian metric
on $G/K$ (compare
Lemma~\ref{lem:transitivity-of-nearby-fixed-points}
~\ref{lem:transitivity:compact-adjoint} below).
By \cite[Section~3]{Carlsson-Goldfarb-coherence}
and \cite[Proposition~3.3]{Ji-asymptotic-Novikov-arithmetic},
$G/K$ has finite asymptotic dimension.~
\footnote{These references only consider connected Lie groups.
          If $G_0$ is the component of the identity, then
          $G_0 / G_0 \cap K \cong G / K$ and the general case
          follows.}
Therefore, we must prove the following result.

\begin{proposition}
\label{prop:uniform-fin-contractible-for-lie-groups}
Let $\Gamma$ be a discrete subgroup of a Lie group $G$
with finitely many components.
Let $K$ be a maximal compact subgroup of $G$.
Then the $\Gamma$-space $G/K$ equipped with a
$G$-invariant Riemannian metric is uniformly
$\Fin$-contractible.
\end{proposition}

It is easy to see that $G/K$ is uniformly
contractible, since $G/K$ is contractible
and has a transitive action by isometries.
Thus, in the torsion-free case, the proof of
Proposition~\ref{prop:uniform-fin-contractible-for-lie-groups}
is trivial.
For the general case, the proof
depends on a fixed point result
(Proposition~\ref{prop:fixed-point-close-by}),
which was explained to the authors by Burkhard Wilking.

The following classical facts about Lie groups will be needed.
\begin{numberlist}
\item[\label{central-by-compact-to-vector-by-compact}]
     Let $G$ be a Lie group with finitely many components
     that has a discrete central subgroup $D$ such that
     the quotient $G/D$ is compact.
     Then $G$ is isomorphic to the semidirect product,
     $V \rtimes K$,
     of a vector group~
     \footnote{Lie groups isomorphic to $\IR^n$ are
               called vector groups.}
     , $V$, with a compact group, $K$,
     acting on $V$.
     This can, for example, be extracted from
     the proof of Lemma XV.3.3 in
     \cite{Hochschild-Lie-groups-book}
     (see the end of the first paragraph on p. 183).
\item[\label{vecter-with-compact-index-is-semidirect}]
     Let $G$ be a Lie group, and let $V$ be a normal
     vector subgroup of $G$ such that $G/V$ is compact.
     Then $G$ is isomorphic to the semidirect product
     $V \rtimes G/V$.
     This follows from
     \cite[Theorem~III.2.3]{Hochschild-Lie-groups-book}.
\item[\label{nonpositively-curved}]
     Let $G$ be a  semisimple Lie group with a finite
     number of components.
     Let $K$ be a maximal compact subgroup of $\Ad_G(G)$,
     and let $L$ be the preimage of $K$ in $G$.
     Then any $G$-invariant metric on $G/L$
     has nonpositive sectional curvature
     (see \cite{Helgason-Book-1962}
     ~\footnote{This is well known.
     Since we did not find the statement in this form,
     we give a proof with precise references:
     We may assume that $G$ is connected,
     since $G_0 / G_0 \cap L \cong G / L$.
     Let $\liel \oplus V$ be a Cartan decomposition
     of the Lie algebra $\lieg$ of $G$ (see
     \cite[III.\S7]{Helgason-Book-1962}).
     Then $\liel$ is a maximal compactly imbedded
     subalgebra of $\lieg$
     (see \cite[Proposition~III.7.4]{Helgason-Book-1962}).
     This means that the subgroup of $\Ad_G( G)$
     corresponding to
     $\liel$ is a maximal compact subgroup.
     Since all maximal compact subgroups are conjugated
     \cite[Theorem~XV.3.1.]{Hochschild-Lie-groups-book},
     we may assume that this subgroup is $K$.
     Thus, the pair $(G,L)$ is of noncompact type
     \cite[p.194/195]{Helgason-Book-1962},
     and any $G$-invariant metric on $G/L$ has
     nonpositive curvature.
     }.)
\end{numberlist}

We make the following convenient definition.

\begin{definition}
\label{def:nearby-fixed-points}
Let $H$ be a closed subgroup of a Lie group $G$.
Assume that there is a $G$-invariant metric on $G/H$.
The pair $(G,H)$
has \emph{nearby $\Fin$-fixed points} if the following holds:

Let $F$ be a finite subgroup of $G$ and $R>0$.
Then there exists a $T>0$ such that for every $F$-invariant
subset $B$ of $G/H$ whose diameter is bounded by $R$,
there is a point fixed by $F$ such that the ball
of radius $T$ around this fixed point contains $B$.
\end{definition}

It is not difficult to check that the existence of $T$,
in the above definition,
does not depend on the chosen metric, provided the metric
is $G$-invariant.

\begin{proposition}
\label{prop:fixed-point-close-by}
Let $G$ be a Lie group with finitely many components and
$K$ a maximal compact subgroup of $G$.
Then $(G,K)$ has nearby $\Fin$-fixed points.
\end{proposition}

The proof of Proposition~\ref{prop:fixed-point-close-by} will
require some preparations.

A Riemannian submersion $f \colon M \to N$ between
Riemannian manifolds is a differentiable surjective map
such that for every $x \in M$, the restriction of
$df \colon T_x M \to T_{f(x)} N$ to
the orthogonal complement of its kernel is an isometry.
Important for us is the following property:
\setlength\labelwidth\origlabelwidth
\begin{numberlist}
\item[\label{lift-for-riemannian-submersion}]
         If $\omega$ is a smooth path in $N$ and $x$ is a lift
         of its the initial point to $M$, then there is a
         canonical
         lift of $\omega$ to $M$ that begins at $x$ and has
         the same arc length as $\omega$.
\end{numberlist}

\begin{lemma}
\label{lem:transitivity-of-nearby-fixed-points}
Let $K \subset H$ be closed subgroups of the Lie group $G$.
\begin{enumerate}
\item \label{lem:transitivity:submersion}
      Suppose that there are $G$-invariant metrics
      on $G/K$ and $G/H$ such that the projection
      $G/K \to G/H$ is a Riemannian submersion.
      If $(H,K)$ and $(G,H)$ have nearby $\Fin$-fixed points,
      then so does $(G,K)$.
\item \label{lem:transitivity:compact-adjoint}
      If $\Ad_G(H)$ is compact, then there are
      $G$-invariant metrics on $G/K$ and $G/H$
      such that the projection
      $G/K \to G/H$ is a Riemannian submersion.
\end{enumerate}
\end{lemma}

\begin{proof}
\ref{lem:transitivity:submersion}
Let $F$ be a finite subgroup of $G$, and
let $B$ be an $F$-invariant subset of $G/K$ of
diameter bounded by $R$.
Since $(G,H)$ has nearby $\Fin$-fixed points
and $\pi$ is non-expanding,
there is a fixed point $gH$ in $G/H$ such
that $\pi(B)$ is contained in ${gH}^T$, the ball of radius
$T$ around $gH$.
Here, $T$ only depends on $F$ and $R$ (and not on $B$).
For $b$ in $B$, choose a geodesic, $\omega_b$, in
$G/H$ from $\pi(b)$ to $gH$.
Let $\varphi(b)$ be the endpoint of the lift of $\omega_b$
to $G/K$ with initial point $b$.
Let $C = \{ \varphi(b) \; | \; b \in B \} \subset \pi^{-1}(gH)$.
Property~\eqref{lift-for-riemannian-submersion}
implies $B \subset C^T$ and $C \subset B^T$.
It follows that the diameter of $C$ with respect to
the metric on $G/H$ is bounded by $R + 2T$.
The diameter with respect to the restriction of the Riemannian
metric to $\pi^{-1}(gH)$ may be larger,
but will still be bounded
by some number $R'$ depending only on $R + 2T$.
This is true because
$H^g = g^{-1} H g$ acts transitively and isometrically with
respect to either distance on $\pi^{-1}(gH)$.
Since $(H,K)$ has nearby $\Fin$-fixed points, there
is a fixed point, $hK$, for $F^g = g^{-1} F g$, and a
$T' > 0$, depending only on $R'$ and $F^g$, such that
$g^{-1}(C) \subset \pi^{-1}(eH) = H/K$ is contained
in $hK^{T'}$.
Therefore, $ghK$ is a fixed point for $F$, and $B$ is contained
in $ghK^{T + T'}$.

\ref{lem:transitivity:compact-adjoint}
Since $\Ad_G(H)$ is compact, there is an
$H$-invariant inner product on the Lie algebra
$\lieg$ of $G$.
The projection $G \to G/H$ induces an
isomorphism from $\lieh^\perp$ to the
tangent space $T_{eH} (G/H)$
of $G/H$ at $eH$, where $\lieh$ is the Lie subalgebra of $\lieg$ corresponding to $H$.
We use this isomorphism to transport the inner product
from $\lieh^\perp$ to $T_{eH} (G/H)$.
Since the inner product on $\lieh^\perp$ is
$H$-invariant, this inner product extends to a
$G$-invariant Riemannian metric on $G/H$.
(Here, it is important that the inner product
on $\lieg$ is $H$-invariant, otherwise the
inner product on $T_{gH} (G/H)$ would depend on the choice
of $g$, and not just on $gH$.)
Similarly, we obtain a $G$-invariant metric on $G/K$.
Thus, the projection, $\pi \colon G/K \to G/H$, is a
Riemannian submersion.
\end{proof}

\begin{proof}
[Proof of Proposition~\ref{prop:fixed-point-close-by}]
We begin by considering two special cases.

Case 1: If $G$ is the semidirect product, $V \rtimes K$,
of a vector group $V$ with a compact group $K$,
then $V \cong G/K$. That is, $G/K$ is Euclidean space
and the result follows from
the Cartan Fixed Point Theorem.

Case 2: Suppose that $G$ semisimple.
Let $L$ be the preimage under $\Ad_G$ of
a maximal compact subgroup of $\Ad_G(G)$
that contains $\Ad_G(K)$.
By Lemma~\ref{lem:transitivity-of-nearby-fixed-points}
~\ref{lem:transitivity:compact-adjoint},
there exists a $G$-invariant metric
on $G/L$, which is nonpositively curved
by \eqref{nonpositively-curved}.
By the Cartan Fixed Point Theorem, $(G,L)$ has nearby
$\Fin$-fixed points.
Since $G/L$ is connected, $L$ will only have a finite
number of components.
Since $G$ is semisimple, the kernel of
$\Ad_G$ is a discrete subgroup
and the intersection of the kernel of $\Ad_G$
with the identity component
$G_0$ of $G$,
which we call $D$, is the center of $G_0$~
\cite[Corollaries~II.5.2+II.6.2]{Helgason-Book-1962}.
Since $L \cap G_0 / D$ is compact,
\eqref{central-by-compact-to-vector-by-compact}
implies that
$L \cap G_0$ is a semidirect product of a vector group
group with a compact group.
Therefore $L \cap G_0 / K \cap G_0 \cong L / K$
is Euclidean space.
Thus, $(L,K)$ has nearby $\Fin$-fixed points
by the Cartan Fixed Point Theorem.
By Lemma~\ref{lem:transitivity-of-nearby-fixed-points},
$(G,K)$ has nearby $\Fin$-fixed points.

For the general case, we proceed by induction on $\dim G$.
If $G$ is not semisimple, there is a closed nontrivial
normal subgroup that is either a vector group $V$
or a torus $T$~
\cite[Lemma~XV.3.6]{Hochschild-Lie-groups-book}.
If the subgroup in question is a torus,
then we may assume that it is contained in $K$. Thus,
the result for $G/T$ implies it for $G$.
If it is a vector group $V$, then consider the
subgroup $VK$ of $G$.
By \eqref{vecter-with-compact-index-is-semidirect},
the first case applies to $VK$.
Of course, $G / VK \cong (G/V) / (VK/V)$.
By induction, $(VK,K)$
and $(G,VK)$ have nearby $\Fin$-fixed points.
Choose a $K$-invariant inner product on
the Lie algebra $\lieg$ of $G$.
Denote by $\liev$ the ideal of $\lieg$ corresponding to $V$.
The projection $G \to G/V$ induces an
isomorphism of the Lie algebra of $G/V$ with $\liev^\perp$
as a $K$-module.
As in the proof of
Lemma~\ref{lem:transitivity-of-nearby-fixed-points}
~\ref{lem:transitivity:compact-adjoint},
we can use left translation to obtain
$G$-invariant, resp.\ $G/V$-invariant, metrics on $G/K$, resp.\
$(G/V)/(VK/K)$, such that the projection
$G/K \to (G/V) / (VK/K)$ is a Riemannian submersion.
Using $G / VK \cong (G/V) / (VK/V)$ and
Lemma~\ref{lem:transitivity-of-nearby-fixed-points}
~\ref{lem:transitivity:submersion},
it follows that $(G/K)$ has nearby $\Fin$-fixed points.
\end{proof}

\begin{lemma}
\label{lem:contract-balls-around-eK}
Let $G$ be a Lie group with finitely many
components and $K$ a maximal compact subgroup.
Equip $G/K$ with a $G$-invariant Riemannian metric.
Then for every $T > 0$ there is an $S > 0$
such that the ball of radius $T$
around $eK$ is $K$-equivariantly contractible
inside the ball of radius $S$.
\end{lemma}

\begin{proof}
By~\cite[Theorem~XV.3.1]{Hochschild-Lie-groups-book},
there is a finite dimensional vector space $V$
with a linear $K$-action and a $K$-equivariant
homeomorphism $G/K \to V$ sending $eK$ to $0$.
Thus,
$eK^T$ is $K$-equivariantly contractible
in $G/K$.
By compactness, this contraction happens
inside some ball of finite radius.
\end{proof}

\begin{proof}
[Proof of Proposition~
  \ref{prop:uniform-fin-contractible-for-lie-groups}]
Let $F$ be a finite subgroup of $G$.
Let $R > 0$ be given.
By Proposition~\ref{prop:fixed-point-close-by}, there is a $T>0$ such that for every $F$-invariant
subset $B$ of $G/K$ whose diameter is bounded by $R$,
there is a point fixed by $F$ such that the ball
of radius $T$ around this fixed point contains $B$.
By Lemma~\ref{lem:contract-balls-around-eK},
there is an $S>0$ such that $eK^T$ is
$K$-equivariantly contractible
inside the ball of radius $S$.

Let
$B'$ an $F$-invariant subset
of $G/K$ whose diameter is bounded by $R$.
Then there is an $F$-fixed point $gK$ in $G/K$ such that $B'$ is contained in $gK^T$.
Since $gK$ is an $F$-fixed point, $F^g = g^{-1} F g$
is a subgroup of $K$.
Thus, $g^{-1} B'$ is contained in $eK^T$ and is $F^g$-invariant. Therefore,
$g^{-1} B'$ is $F^g$-equivariantly contractible
in $eK^S$.
Applying $g$, this means that $B'$ is $F$-equivariantly
contractible in $gK^S$.
\end{proof}


\section{Open covers and simplicial compexes}
\label{sec:open-covers-simplicial-complexes}

A map $f \colon X \to Y$ between metric spaces is
{\em metrically coarse} if it is proper
and satisfies the following growth condition:
for all $R>0$ there is an $S>0$ such that
\begin{eqnarray*}
  \quad d_X(x,y) < R \quad
   \Longrightarrow \quad d_Y(f(x),f(y)) < S.
\end{eqnarray*}
Two such maps, $f$ and $g$, are said to be {\em bornotopic}
if there is a constant
$C > 0$ such that $d_Y(f(x),g(x)) < C$ for all $x$ in $X$
(see \cite[Section 2]{Higson-Roe-On-coarse-Baum-Connes}).
A metrically coarse  homotopy between proper continuous maps
is called a {\em metric homotopy}.
In particular, a metric homotopy is a bornotopy
and a proper homotopy.

For the following lemma compare
\cite[Lemma 3.3]{Higson-Roe-On-coarse-Baum-Connes}.

\begin{lemma}
\label{lem:extend-coare-continuously}
Let $X$ be a proper metric space that has a finite dimensional
$CW$-structure such that the diameters of its cells are
uniformly bounded, and
let $Y$ be a uniformly contractible proper metric space.
Let $A$ be a subcomplex of $X$ containing the $0$-cells of $X$.
Then every continuous metrically coarse map
$f_0 \colon A \to Y$ can be extended to a
continuous metrically coarse map
$f \colon X \to Y$.
\end{lemma}

\begin{proof}
Let $X_n$ be the union of $A$ with the $n$-skeleton of $X$.
Since $X_n$ is obtained from $X_{n-1}$ by attaching $n$-cells,
we extend $f_0$ inductively to $f_n \colon X_n \to Y$.
For every $n$-cell, $e$, not in $A$, we must extend to all
of $e$ the restriction of $f_{n-1}$
to the boundary, $\dd e$, of $e$.
Because $f_{n-1}$ is metrically coarse and
cells in $X$ have uniformly bounded diameter,
$f_{n-1}(\dd e)$ has uniformly bounded diameter.
Since $Y$ is uniformly contractible, there is a
continuous metrically coarse extension of $f_{n-1}$ to $X_n$.
This finishes the construction of $f$.

It remains to show that $f$ is proper.
Since the diameters of cells in $X$ are uniformly bounded
and $X$ is finite dimensional, there is an $R > 0$ such
that $X = (X^0)^R$, where $X^0$ denotes the $0$-skeleton
of $X$.
Since $X$ is a proper metric space, this implies that a
subset $B$ of $X$ has compact closure if and only if
$B^R \cap X^0$ is finite.
Let $B_\alpha$ be a closed ball of radius $\alpha$ in $Y$.
Because $f$ is metrically coarse, there is an $S > 0$
such that $(f^{-1}(B_\alpha))^R$ is contained in
$f^{-1}((B_\alpha)^S)$.
Since $(B_\alpha)^S$ is compact, $f_0$ is proper,
and $X^0$ is closed in $X$,
it follows that
$(f^{-1}(B_\alpha))^R \cap X^0 \subset
  f^{-1}((B_\alpha)^S) \cap X^0$ is finite.
Since $f^{-1}(B_\alpha)$ is closed, this implies that
$f^{-1}(B_\alpha)$ is compact.
Therefore, $f$ is proper.
\end{proof}

If $\calu$ is an open cover of a space $X$, then the
realization of its nerve, $|\calu|$, is a simplicial complex.
Denote by $[U]$ the vertex of $|\calu|$ corresponding to
$U$ in $\calu$.
A partition of unity subordinate to $\calu$,
$(\varphi_U)_{U \in \calu}$, induces a map
$g \colon X \to | \calu |$ defined by:
\begin{equation}\label{induced}
g (x) = \sum_{U \in \calu} \frac{\varphi_U(x)[U]}
                          {\sum_{V \in \calu} \varphi_V(x)}
\end{equation}

The {\it Euclidean path length metric}
~\footnote{In \cite{Yu-finite-asymptotic-dimension-Novikov}
        and~\cite{Bartels-finite-asymp-dim}
           the spherical metric has been used.
           Using the Euclidean metric is convenient for
           computation in
           Proposition~\ref{lem:estimates-for-map-to-nerve}.
           In any event, the difference between the
           Euclidean and
           the spherical metric is not important for
           this paper.}
on a simplicial complex
is the unique path length metric that restricts to the standard
Euclidean metric on each simplex.
For an open cover $\calu$, we will always equip $|\calu|$
with the Euclidean path length metric.

\begin{proposition}
\label{prop:bornotopic}
Let $X$ be a complete proper path metric space and
let $\calu$ be a locally finite, bounded, finite dimensional
cover of $X$
whose Lebesgue number is positive.
If $g \colon X \to |\calu|$ is induced
by a partition of unity subordinate to $\calu$
as in \eqref{induced},
then $g$ is a bornotopy equivalence.
\end{proposition}

\begin{proof}
This follows from
\cite[Section 3]{Roe-hyperbolic-spaces-exotic-cohomology}
(see also
\cite[Proposition~3.2]{Higson-Roe-On-coarse-Baum-Connes}).
In these references the spherical metric rather than the
Euclidean metric is used, but since $|\calu|$
is finite dimensional,
this distinction is not important.
\end{proof}

\begin{lemma}
\label{lem:right-inverse-for-g-calu}
Let $X$ be a uniformly contractible complete proper path metric space.
Assume that $X$ has the structure of a finite
dimensional $CW$-complex.
Let $\calu$ be a locally finite, bounded,
finite dimensional cover of $X$
whose Lebesgue number is positive,
and let $g \colon X \to |\calu|$ be induced by a
partition of unity subordinate to $\calu$ as in \eqref{induced}.
Then $g$ has a right homotopy inverse up to
metric homotopy. That is, there is a continuous
and metrically coarse map $f \colon |\calu| \to X$
and a metric homotopy $H \colon X \x [0,1] \to X$
from $f \circ g$ to $\id_X$.
\end{lemma}

\begin{proof}
By refining the $CW$-structure if necessary, we can assume that
the cells in $X$ have uniformly bounded diameter.
By Proposition~\ref{prop:bornotopic},
$g$ has a bornotopy inverse $f$.
Using Lemma~\ref{lem:extend-coare-continuously},
we can assume that $f$ is also continuous.
The existence of $H$ follows by applying
Lemma~\ref{lem:extend-coare-continuously}
to the subspace $X \x \{ 0, 1 \}\subset X \x [0,1]$
and the map that is $f \circ g$
on $X \x \{ 0 \}$ and the identity of $X$
on $X \x \{ 1 \}$.
\end{proof}

Although the Euclidean path length metric changes
under restriction
to subcomplexes, there is the following estimate.

\begin{lemma}
\label{lem:compare-l^2-and-d-pl}
Let $\sigma$ and $\tau$ be intersecting
faces of the $k$-simplex $\Delta$.
Denote by $d_\pl$ the Euclidean path length metric
on the subcomplex of $\Delta$ spanned by $\sigma$ and $\tau$,
and denote by $d_{\Delta}$ the Euclidean standard metric on
$\Delta$. Then for $x$ in $\sigma$ and $y$ in $\tau$
\[
d_\pl ( x,y) \leq 3 \sqrt{k+1} \cdot d_{\Delta} (x,y).
\]
\end{lemma}

\begin{proof}
Let $v_0$ be a vertex that is contained in $\sigma \cap \tau$.
For a face $\rho$ of $\Delta$ that contains $v_0$, define the
projection $p_\rho \colon \Delta \to \rho$ as the simplicial map
which is the identity on $\rho$ and maps all vertices not
in $\rho$ to $v_0$.
Using the Cauchy-Schwarz inequality,
\[
d_\Delta ( p_\rho(z), p_\rho(z') ) \leq
\sqrt{k+1} \cdot d_\Delta ( z, z')
               \quad \forall z,z' \in \Delta.
\]
Note that $p_\sigma(x) = x$, $p_\tau(y) = y$,
$p_{\sigma \cap \tau} (x) = p_\tau (x)$
and $p_{\sigma \cap \tau} (y) = p_\sigma (y)$.
Therefore,
\begin{eqnarray*}
d_\pl (x,y)  & \leq  & d_\pl ( x, p_{\sigma \cap \tau} (y) ) +
                       d_\pl ( p_{\sigma \cap \tau} (y),
                               p_{\sigma \cap \tau} (x) )    +
                       d_\pl ( p_{\sigma \cap \tau} (x), y ) \\
             & =     & d_\pl ( p_{\sigma} (x), p_\sigma (y) ) +
                       d_\pl ( p_{\sigma \cap \tau} (y),
                               p_{\sigma \cap \tau} (x) )    +
                       d_\pl ( p_{\tau} (x), p_{\tau} (y) ) \\
             & =     & d_\Delta ( p_{\sigma} (x), p_\sigma (y) ) +
                       d_\Delta ( p_{\sigma \cap \tau} (y),
                                  p_{\sigma \cap \tau} (x) ) +
                       d_\Delta ( p_{\tau} (x), p_{\tau} (y) ) \\
             & \leq  & 3 \sqrt{k} \cdot d_\Delta ( x, y).
\end{eqnarray*}
\end{proof}

The finite asymptotic dimension of a metric space is
equivalent to the existence of certain contracting maps
to finite dimensional simplicial complexes
(see \cite[p.30]{Gromov-asymptotic-invariants}).
We will use the following version of this (compare
\cite[Lemma 6.3]{Yu-finite-asymptotic-dimension-Novikov}).

\begin{lemma}
\label{lem:estimates-for-map-to-nerve}
Let $\calu$ be an $n$-dimensional bounded cover
of a path length metric space $X$.
Assume that the Lebesgue number $R$ of $\calu$ is positive.
Then there is a partition of unity subordinate to $\calu$
such that for the induced map
$g \colon X \to |\calu|$~\eqref{induced},
\[
d_{\pl}( g(x),g(y) ) \leq C_n \frac{d(x,y)}{R}
\]
for all $x$, $y \in X$.
Here, $d_\pl$ denotes the Euclidean path length metric
on $|\calu|$, and $C_n$ is a constant that depends
only on $n$.
\end{lemma}

\begin{proof}
Let $V_\calu$ be the vector space of sequences of real numbers
indexed by $\calu$.
There is a canonical embedding $|\calu| \to V_\calu$.
We will use both the $l^1$-norm, $\onenorm{\cdot}$,
and the $l^2$-norm, $\twonorm{\cdot}$, on $V_\calu$.
For $x$ in $X$ and $U$ in $\calu$, let $x_U = d(x,X-U)$
and $f(x) = (x_U)_{U \in \calu} \in V_\calu$.
Define a partition of unity subordinate to $\calu$
by $\varphi_U(x) = \frac{x_U}{ \onenorm{f(x)} }$.
The map $g \colon X \to |\calu| \subset V_\calu$
induced by this partition of
unity is given by
$x \mapsto \frac{f(x)}{\onenorm{f(x)}}$.

Let $x$, $y \in X$ be given.
Clearly, $\twonorm{f(x)} \leq \onenorm{f(x)}$.
Because $R$ is the Lebesgue number of $\calu$, there
is at least one $U$ in $\calu$ for which $y_U \geq R$.
Therefore, $\onenorm{f(y)} \geq R$.
Since $\calu$ is $n$-dimensional, there are at most $2(n+1)$
members of $\calu$ for which $x_U \neq 0$ or $y_U \neq 0$.
{}From $|x_U - y_U| \leq d(x,y)$, we conclude
$\big| \onenorm{f(y)} - \onenorm{f(x)} \big|
                           \leq 2(n+1) \, d(x,y)$
and
$\twonorm{ f(x) - f(y) } \leq \sqrt{2(n+1)}\, d(x,y)
                         \leq 2(n+1)\, d(x,y)$.
Using these estimates, it follows that
\begin{align*}
\twonorm{g(x) - g(y)}
   & = \twonorm{ \frac{f(x)}{\onenorm{f(x)}} -
                       \frac{f(y)}{\onenorm{f(y)}} }
   \\
   & \leq \twonorm{ \frac{ \onenorm{f(y)} f(x) -
                                      \onenorm{f(x)} f(x)}
                         { \onenorm{f(x)} \onenorm{f(y)} } }
       + \twonorm{ \frac{ \onenorm{f(x)} f(x) -
                                       \onenorm{f(x)} f(y)}
                         { \onenorm{f(x)} \onenorm{f(y)} } }
   \\
   & \leq \frac{ \big| \onenorm{f(y)} - \onenorm{f(x)} \big|
                                          \twonorm{f(x)} }
               { \onenorm{f(x)} \onenorm{f(y)} }
       + \frac{ \twonorm{ f(x) - f(y) } }
               { \onenorm{ f(y) } }
   \\
   & \leq \frac{4(n+1) d(x,y) }{R}.
\end{align*}
Restricted to simplices, the embedding $|\calu| \to V_\calu$
is an isometry (using the $l^2$-norm on $V_\calu$).
If $d(x,y) < R$, then there is a
$U$ in $\calu$ containing $x$ and $y$,
since $R$ is the Lebesgue number of $\calu$.
Thus, there are simplices $\sigma$ and $\tau$ of $|\calu|$
such that $g(x) \in \sigma$,
$g(y) \in \tau$, and $[U]\in\sigma \cap \tau$.
In particular, $\sigma \cap \tau$ is not empty.
Because $\sigma$ and $\tau$ have dimension less or equal to $n$,
they span a  simplex $\Delta$ in $V_\calu$
whose dimension is at most $2n$.
Since the $l^2$-norm gives the Euclidean metric on $\Delta$,
we can use Lemma~\ref{lem:compare-l^2-and-d-pl}
to estimate
\[
d_\pl(g(x),g(y)) \leq 3 \sqrt{2n+1} \cdot \twonorm{g(x) - g(y)}
      \leq 12 \sqrt{2n + 1} (n+1)\frac{d(x,y)}{R},
\]
whenever $d(x,y) \leq R$.
Since the metric on $X$ is a path length metric,
the inequality follows for all $x$, $y \in X$.
\end{proof}


\section{Controlled algebra}\label{sec:controlled-algebra}

Let $X$ be a proper metric space and $\cala$ a small
additive category.
Let $\Gamma$ be a group that acts on $X$ by isometries
and on $\cala$ by additive functors.
(We will consider only left actions.)
Fundamental to this paper will be the additive category
$\cala(X)$ of certain continuously controlled modules.
This is a minor variation of similar constructions from the
literature, the first of which appeared
in~\cite{ACFP-continuous-control}.

Let $Z = X \x \Gamma \x [0,1)$.
An object, $M$, in $\cala(X)$ is given by a sequence of
objects $(M_z)_{z \in Z}$ in $\cala$,
subject to the conditions:
\begin{numberlist}
\item The image of
      $\supp M = \big\{ z \, \big| \, M_z \neq 0 \big\}$
      under the projection $Z \to X \x [0,1)$ is locally finite.
\item For every $x \in X$ and $t \in [0,1)$,
$\supp M \cap \big(\{ x \} \x \Gamma \x \{t\}\big)$
      is finite.
\end{numberlist}

A morphism, $\varphi \colon M \to N$, in $\cala(X)$ is
given by a sequence of morphisms,
$(\varphi_{z,z'} \colon M_{z'} \to M_{z})_{(z,z') \in Z^2}$,
in $\cala$, subject to the conditions:
\begin{numberlist}
\item   $\supp \varphi = \big\{ (z,z')\, \big|\,
                     \varphi_{z,z'} \neq 0 \big\}$
        is {\it continuously controlled} at $X \x \{ 1 \}$.
        That is, for every $x \in X$ and every open neighborhood
        $U$ of $(x,1)$ in $X \x [0,1]$,
        there is a (smaller) open
        neighborhood $V$ of $(x,1)$ in $X \x [0,1]$ such that
        $(X \x [0,1) - U) \x V$ and
        $V \x (X \x [0,1) - U)$ do not intersect the image
        of $\supp \varphi$ under the projection
        $Z \to X \x [0,1]$.
\item   For a fixed $z$ in $Z$,
        $\big\{ z' \,\big|\, (z,z')\in \supp \varphi
            \;\, {\rm or} \;\, (z',z)\in\supp \varphi \big\}$
        is finite.
\item[\label{metric-condition-for-morphisms}]
       There exists an $R > 0$ such that
       $((x,\gamma,t),(x',\gamma',t')) \in \supp \varphi$
       implies $d(x,x') < R$.
\end{numberlist}
Composition of morphisms is given by the usual matrix
multiplication.

If $\calx$ is a collection of subsets of $Z$, then we
define $\cala(\calx)$ as the full subcategory of $\cala(X)$
whose objects, $M$, satisfy the additional condition:
\begin{numberlist}
\item There is an $S \in \calx$ such that $\supp M \subset S$.
\end{numberlist}
If $\calx$ is closed under finite unions, then $\cala(\calx)$
is again an additive category.
If $\caly$ is another collection of subsets of $Z$ such that
for every $S \in \caly$ there is a $T \in \calx$ such that
$S \subset T$, then $\cala(\caly)$ is a subcategory of
$\cala(\calx)$.
Furthermore, $\cala(\caly)$ will define a Karoubi filtration
\cite[Definition 1.27]
{Carlsson-Pedersen-Controlled-algebra-Novikov}
of $\cala(\calx)$ if, in addition, the following is satisfied:
\begin{numberlist}
\item For every $S \in \caly$ and morphism $\varphi$
      in $\cala(\calx)$, there is a $T \in \caly$ such that
      $S^\varphi = \big\{ z \; \big| \; \exists \, z' \in S \;
          \mbox{such that} \; \varphi_{z,z'} \neq 0 \big\}$
      is contained in $T$.
\end{numberlist}
The quotient of this Karoubi filtration will be denoted by
$\cala(\calx, \caly)$.

Clearly, $\cala(X) = \cala( \{ X \})$
and $\cala(\calx) = \cala( \calx, \emptyset)$.
If $\calx$ and $\caly$ are $\Gamma$-invariant, then
the formula $(g(M))_z = g(M_{g^{-1}z})$
defines an action of $\Gamma$ on
$\cala(\calx, \caly)$.
For a subgroup $G$ of $\Gamma$, denote
the corresponding fixed point category by
$\cala^G(\calx,\caly)$.
It is not difficult to check that taking Karoubi quotients
and taking fixed categories commute.
Therefore, $\cala^G( \calx, \caly)$ is the quotient
of $\cala^G( \calx)$ by $\cala^G(\caly)$.

Let $p \colon X \to X'$ be a continuous map.
For a closed subset $Y$ of $X$, let
$\calx(Y,p)$ be the collection of subsets $S$ of $Z$
with the following properties:
\begin{numberlist}
\item [\label{define-calx(Y,p)-limit}]
      If $(x,1)$ is a limit point of the image of
      $S$ under the projection $Z \to X \x [0,1)$,
      then $x \in Y$.
\item [\label{define-calx(Y,p)-bounded}]
      There is an $R > 0$ such that the image of $S$
      under the projection $Z \to X$ is contained in $Y^R$.
\item [\label{define-calx(Y,p)-compact}]
      There is a compact set $K_0 \subset X'$ such that
      the image of $S$ under the composition of the
      projection $Z \to X$ with $p$ is contained in $K_0$.
\end{numberlist}
Let $\calx(Y,p)_0$ be the collection of subsets $S$ of $Z$
that satisfy \eqref{define-calx(Y,p)-bounded},
\eqref{define-calx(Y,p)-compact} and
the following strengthening of \eqref{define-calx(Y,p)-limit}:
\begin{numberlist}
\item The set of limit points of the image of
      $S$ under the projection $Z \to X \x [0,1)$
      is disjoint from $X \x \{ 1 \}$.
\end{numberlist}

Let $G$ be a subgroup of $\Gamma$ such that $Y$
is $G$-invariant.
The controlled categories that we will use are:
\begin{align*}
\cala_p^G(Y)          & = \cala^G(\calx(Y,p)), \\
\cala_p^G(Y)_0        & = \cala^G(\calx(Y,p)_0), \\
\cala_p^G(Y)^\infty   & = \cala^G(\calx(Y,p), \calx(Y,p)_0), \\
\cala_p^G(X,Y)        & = \cala^G(\calx(X,p), \calx(Y,p)).
\end{align*}
Note that the definitions of these categories depend
on the metric space $X$, even though this is not reflected
in the
notation.
However, because $X$ is a proper metric space,
changing $X$ to a smaller or larger proper metric space
that still contains $Y$ will only change the categories
up to $G$-equivariant equivalence.
Since all of our metric spaces will be proper metric spaces,
we can disregard this dependence on $X$.
It should also be noted that
condition~\eqref{metric-condition-for-morphisms}
is only important in the definitions of $\cala_p^G(Y)$ and
$\cala_p^G(Y)_0$ and not in the definition of
$\cala_p^G(Y)^\infty$,
where it affects the category only
up to equivalence (compare
\cite[Lemma~3.15]{Bartels-finite-asymp-dim}).
Similarly, our notation is slightly imprecise because
even before taking fixed categories these
categories depend on the group we have in mind.
But this dependence also only changes the category
up to equivariant equivalence and can safely be ignored.

Let $\mathbb{K}^{-\infty}$ denote the functor from the
category of small additive categories to the
category of spectra, which assigns an additive category to
its associated non-connective $K$-theory
spectrum~\cite{Pedersen-Weibel-delooping}.
A crucial fact is that applying $\mathbb{K}^{-\infty}$
to a Karoubi filtration produces a fibration of spectra,
which induces a long exact
sequence in $K$-theory
\cite[Theorem~1.28]{Carlsson-Pedersen-Controlled-algebra-Novikov}.
By definition, the two following sequences
are Karoubi sequences.
That is, the final category
is a Karoubi quotient of the middle category by the first.
Therefore, each induces a long exact sequence in $K$-theory.
\begin{numberlist}
\item [\label{germs-at-infty-sequence}]
      $\cala^G_p(X)_0 \to \cala^G_p(X) \to \cala^G_p(X)^\infty$
\item [\label{pair-sequence-for-cala}]
      $\cala^G_p(Y) \to \cala^G_p(X) \to \cala^G_p(X,Y)$
\end{numberlist}
If, in addition, $Y^R = X$ for some $R > 0$, then
\begin{numberlist}
\item [\label{pair-sequence-for-cala^infty}]
      $\cala^G_p(Y)^\infty \to \cala^G_p(X)^\infty
                           \to \cala^G_p(X,Y)$
\end{numberlist}
is also Karoubi sequence.
To see this, consider the following commutative diagram
in which the first two columns and the second two rows are
Karoubi filtrations.
\[
\xymatrix{
    \cala^G_p(Y)_0
       \ar[d] \ar[r]^{\cong} &
    \cala^G_p(X)_0
       \ar[d] &
    \\
    \cala^G_p(Y)
       \ar[d]_a \ar[r] &
    \cala^G_p(X)
       \ar[d]_b \ar[r] &
    \cala^G_p(X,Y)
    \ar[d]_c
    \\
    \cala^G_p(Y)^\infty
      \ar[r] &
    \cala^G_p(X)^\infty
      \ar[r] &
    Q
}
\]
The assumption $Y^R = X$ implies that
$\cala^G_p(Y)_0 \to \cala^G_p(X)_0$ is an equivalence
of categories.
A short exercise in the definition of Karoubi
filtrations shows that $c$ is an equivalence of categories.

In \cite{Weiss-excision-and-restriction},
it is proven that $X \mapsto K_*(\cala_{\id_X} (X)^\infty)$
is a locally finite homology theory.
In particular, it is homotopy invariant.
It follows from the definition that $\cala_{\id_X}(X)_0$
is functorial in $X$ for metrically coarse maps and that
bornotopic maps induce the same map in $K$-theory.
These facts imply:
\begin{numberlist}
\item [\label{homotopy-invariance-without-germs}]
      The category $\cala(X)$ is functorial in $X$ for
      metrically coarse continuous maps
      (compare~\cite[Remarks~3.5, 3.10]{Bartels-finite-asymp-dim}).
      If $f,g \colon X \to Y$ are two such maps
      that are metrically homotopic
      (see Section~\ref{sec:open-covers-simplicial-complexes}),
      then $f$ and $g$ induce the same map from
      $K_* (\cala(X))$ to $K_* (\cala(Y))$.
\end{numberlist}

In \cite[Definition~2.1]{BR-coefficients},
 the additive category $\cala \ast_\Gamma T$
is constructed for a
$\Gamma$-set $T$.
By \cite[Proposition~2.8(iii)]{BR-coefficients},
there is a canonical equivalence of categories
$\cala \ast_H H/H \to \cala \ast_G G/H$,
for every subgroup $H$ of $\Gamma$.
If $\cala$ is the category of finitely generated free
$R$-modules for a ring $R$ and the action of $H$ on
$\cala$ is trivial, then $\cala \ast_H H/H$ is equivalent
to the category of finitely generated free $R[H]$-modules.
For this reason, we will denote the
category $\cala \ast_H H/H$ by $\cala[H]$
for an arbitrary $\cala$, and
call the objects of this category $\cala[H]$-modules.
In \cite[Section~3]{BR-coefficients},
the $\Or (\Gamma)$-spectrum
$\bfK_\cala$ is also introduced
(using the $\mathbb{K}^{-\infty}$-functor).
Associated to it is the assembly map
\begin{equation}
\label{assembly-map-with-coeficients}
H_*(\Ebar; \bfK_\cala) \to H_*(\pt;\bfK_\cala)
                       = K_*(\cala[\Gamma]).
\end{equation}
If $\cala$ is the category of finitely generated free
$R$-modules, then this is the assembly map
\eqref{K-assembly-map-ring}.
Let $p_\Gamma \colon \Ebar \to \Ebar/\Gamma$ denote the
quotient map.
It is not hard to check that $\cala[\Gamma]$ is equivalent
to the category $\cala_{p_\Gamma}^\Gamma(\Ebar)_0$.
Controlled algebra can be used to describe assembly maps
as forget control maps.
An instance of this is \cite{Hambleton-Pedersen-identify},
where the continuously controlled
category $\calb_{\Gamma}(\Ebar \x [0,1);R)$ is used to identify
various versions of the assembly map.
If one modifies the definition of this category to
allow for coefficients in $\cala$, then, since $\Gamma$
acts on $\Ebar$ with finite isotropy,
the fixed point category
$\cala^\Gamma_{p_\Gamma}(\Ebar)^{\infty}$ can be identified with
the continuously controlled category
$\calb_{\Gamma}(\Ebar \x [0,1);\cala)^{>0}$.
The only difference between the two categories is
condition~\eqref{metric-condition-for-morphisms},
but as mentioned before,
this changes the category only up to equivalence.
Therefore,
\cite[Theorem 7.4]{Hambleton-Pedersen-identify}
implies the following fact.
\begin{numberlist}
\item [\label{assembly-as-forget-control}]
      If $\Ebar$ is equipped with
      a $\Gamma$-invariant metric,
      then the boundary map in the long exact sequence
      associated to \eqref{germs-at-infty-sequence}, with
      $X = \Ebar$,
      is equivalent to the assembly map
      \eqref{assembly-map-with-coeficients}.
\end{numberlist}


\section{A vanishing result}\label{sec:vanishing}
In this section, a key component of the proof
of Theorem~\ref{thm:split-injectivity-with-coeficients}
is established, namely Proposition~\ref{prop:vanish} below.
\begin{proposition} \label{prop:squeezing}
Let $X$ be a proper metric space.
For each $n \in\IN$, let $Q_n$ be a  simplicial complex
equipped with the Euclidean path length metric and
$g_n \colon X \to Q_n$ a continuous and metrically coarse map
satisfying the following:
For every $R > 0$ there exists an $S > 0$ such that
for all $x$ and $y$ in $X$ with $d(x,y) < R$,
$d(g_n(x),g_n(y)) < \frac{S}{n}$.
Then for any $a \in K_m(\cala(X))$
there is an $n_0 = n_0(a)$ such that
$(g_n)_*(a) = 0 \in K_m(\cala(Q_n))$ for all $n \geq n_0$.
\end{proposition}

\begin{proof}
This is almost \cite[Corollary 4.3]{Bartels-finite-asymp-dim}.
In \cite{Bartels-finite-asymp-dim}, the spherical metric
rather then the Euclidean metric is used,
but this does not affect the argument.
It is only important that the metric is the same on
every simplex.
\end{proof}

\begin{proposition}
\label{prop:vanish}
Let $X$ be a uniformly contractible,
complete, proper, path length metric space of finite asymptotic
dimension.
Assume that $X$ has the structure of a finite dimensional
$CW$-complex.
Let $\calx$ be a collection of closed subsets, $S$,
of $X\x [0,1)$, where $S=K \x [0,1)$ for some closed
subset $K$ of $X$, that is closed under finite unions,
and assume that for every closed subset $K$ of $X$
and every $\alpha > 0$, $K^\alpha \x [0,1)$
is contained in $\calx$.
Then the $K$-theory of $\cala(\calx)$ vanishes.
\end{proposition}

\begin{proof}
Since $X$ has finite asymptotic dimension,
there is a sequence of bounded open covers $\calu_n$ of $X$
such that the Lebesgue number of $\calu_n$ exceeds $n$.
By Lemma~\ref{lem:finite-asymp-dim-and-locally-finite},
we can assume that each $\calu_n$ is locally finite.
By Lemma~\ref{lem:estimates-for-map-to-nerve} and
Proposition~\ref{prop:squeezing}, there is a sequence of maps,
$g_n \colon X \to |\calu_n|$, induced by partitions of unity,
such that for every closed subset $Y$ of $X$ and every
$a \in K_*(\cala(Y))$, there is an $n_0$ such that
$(g_n)_*(a) = 0$ for $n \geq n_0$.
By Lemma~\ref{lem:right-inverse-for-g-calu},
each $g_n$ is invertible up to metric homotopy.
Thus, for each $n$, there is a
continuous metrically coarse map
$f_n \colon |\calu_n| \to X$ and a metric homotopy
$H_n \colon X \x [0,1] \to X$ from $f_n \circ g_n$ to $\id_X$.

Let $b \in K_*(\cala(\calx))$ be given.
Choose a closed subset $K$ of $X$ such that $b$ is the
image of some $a \in K_*(\cala(K))$ under the inclusion
$\iota_K \colon \cala(K) \to \cala(\calx)$.
Let $Q_n$ be the smallest subcomplex of $|\calu_n|$
that contains $g_n(K)$.
Since $f_n$ and $H_n$ are metrically coarse, there
is an $\alpha_n > 0$ such that
$f_n(Q_n)$, $H_n(K \x [0,1]) \subset K^{\alpha_n}$.
Consider the restriction maps
$g_n|_K \colon K \to Q_n$ and
$f_n|_{Q_n} \colon Q_n \to K^{\alpha_n}$.
Then the restriction of $H_n$ to $K \x [0,1]$
gives a metric homotopy
from $f_n|_{Q_n} \circ g_n|_{K}$ to the inclusion
$i_{\alpha_n} \colon K \to K^{\alpha_n}$.
{}From \eqref{homotopy-invariance-without-germs}
we conclude that
$(f_n|_{Q_n} \circ g_n|_{K})_* (a) = (i_{\alpha_n})_*(a)$.
By Proposition~\ref{prop:squeezing},
eventually $(g_n|_K)_* (a) = 0$.
The conclusion now follows since $i_{\alpha_n}(a)$
also maps to $b$ under the inclusion
$\iota_{K^{\alpha_n}} \colon
         \cala(K^{\alpha_n}) \to \cala(\calx)$.
\end{proof}


\section{The Descent Principle and Homotopy Fixed Sets}
\label{sec:descent}

Let $\bfS$ be a spectrum with $\Gamma$-action.
For a space $X$ with $\Gamma$-action, let
$\Map_\Gamma(X,\bfS)$ be the spectrum whose
$n$-th space is $\Map_\Gamma(X,\bfS_n)$,
where $\bfS_n$ is the $n$-th space in the spectrum $\bfS$.
The fixed point spectrum, $\bfS^{\Gamma}$,
is defined to be $\Map_{\Gamma}(\pt,\bfS)$.
Let $\calf$ be a family of subgroups of $\Gamma$ that is
closed under conjugation and taking subgroups.
If $X$ is a space with $\Gamma$-action,
then its $\calf$-homotopy fixed point set is defined by
$X^{h_\calf \Gamma} = \Map_\Gamma(E_\calf \Gamma, X)$.
In general, spectra can be difficult to work with,
however, many arguments for spaces can be extended to
$\Omega$-spectra by applying them levelwise.
A particularly useful property they possess is that
the homotopy groups of a product of $\Omega$-spectra
is the product of the homotopy groups.
In order to define the $\calf$-homotopy fixed
points of a spectrum
with $\Gamma$-action, recall that there is a
fibrant replacement functor,
$\bfR:{\rm { SPECTRA} \to \Omega-{ SPECTRA}}$,
that comes equipped with
a natural weak equivalence $\bfA \to \bfR(\bfA)$.
For a construction of this functor, see
for example,~\cite[Section 2]{Lueck-Reich-Varisco-commuting}.
If $\bfS$ is a spectrum with $\Gamma$-action
then so is $\bfR(\bfS)$,
and it is not difficult to check that $\bfR$
commutes with taking
fixed points.
That is, $(\bfR(\bfS))^\Gamma = \bfR(\bfS^\Gamma)$.
The {\it $\calf$-homotopy fixed point spectrum}
\footnote{Both authors of this paper are guilty of
          stating an incorrect
          definition of the homotopy fixed point spectrum
          (forgetting $\bfR$)
          \cite{Bartels-on-the-domain},
          \cite{Rosenthal-splitting-K-theory}.
          However, in those papers, homotopy fixed points are
          only
          applied to $\Omega$-spectra. In this case, both
          definitions yield weakly equivalent spectra.
          Thus, the main results are not
          affected.
          We thank Holger Reich for pointing out
          that the correct
          definition of the homotopy fixed point spectrum
          must involve the functor $\bfR$.}
is defined to be the spectrum
\[
\bfS^{h_\calf \Gamma} =
    \Map_{\Gamma}( E_\calf \Gamma , \bfR(\bfS) ).
\]
Notice that the projection $E_\calf \Gamma \to \pt$ induces a
natural transformation $\bfS^\Gamma \to \bfS^{h_\calf \Gamma}$.

\begin{lemma}
\label{lem:on-homotopy-fixed-spectra}
Let $\calf$ be a family of subgroups of $\Gamma$.
\begin{enumerate}
\item \label{lem:on-h-f-spectra:equivalence}
     If $F \colon \bfS \to \bfT$ is an equivariant
     map of spectra with $\Gamma$-action such that
     $F^G \colon \bfS^G \to \bfT^G$ is a weak
     homotopy equivalence for
     every subgroup $G$ in $\calf$,
     then $\bfS^{h_\calf \Gamma} \simeq \bfT^{h_\calf \Gamma}$.
\item \label{lem:on-h-f-spectra:coinduction}
     Let $\bfB$ be an $\Omega$-spectrum with a $G$-action,
     where $G$ is in $\calf$,
     and let $\Gamma$ act on $\bfS = \Map_G(\Gamma,B)$
     by $(\gamma f)(x)=f({\gamma}^{-1}x)$,
     where $\gamma \in \Gamma$.
     Then $\bfS^{\Gamma} \simeq \bfS^{h_\calf \Gamma}$.
\end{enumerate}
\end{lemma}

\begin{proof}
\ref{lem:on-h-f-spectra:equivalence}
The proof of \cite[Lemma~4.1]{Rosenthal-splitting-K-theory}
shows that the corresponding result holds for spaces,
but this also implies the result for spectra.
For every $G \in \calf$ and $n \in \IN$,
$(\bfR(\bfS)_n)^G \to (\bfR(\bfT)_n)^G$ is a weak equivalence
because $\bfR$ commutes with fixed points.
Thus, $(\bfR(\bfS)_n)^{h_\calf \Gamma} \to
      (\bfR(\bfT)_n)^{h_\calf \Gamma}$
is weak equivalence for all $n \in \IN$.

\ref{lem:on-h-f-spectra:coinduction}
The proof of \cite[Lemma~4.2]{Rosenthal-splitting-K-theory}
shows that the corresponding result holds for spaces.
The statement for spectra can be deduced from this as follows.
A product of $\Omega$-spectra is again an $\Omega$-spectrum.
In particular, $\bfS^H$ is an $\Omega$-spectrum for every
subgroup $H$ of $\Gamma$.
Thus, $(\bfS_n)^H \to (\bfR(\bfS)_n)^H$ is a weak
equivalence for every $n$.
Therefore, $(\bfS_n)^{h_\calf \Gamma} \to
      (\bfR(\bfS)_n)^{h_\calf \Gamma}$
is a weak equivalence for every $n$.
On the other hand,
$(\bfS_n)^\Gamma \to (\bfS_n)^{h_\calf \Gamma}$
is a weak  equivalence
by the space version of \ref{lem:on-h-f-spectra:coinduction}.
\end{proof}

\begin{lemma}
\label{lem:products}
Let $\calf$ be the family of finite subgroups
of $\Gamma$, $H$ a finite subgroup of $\Gamma$,
$\calb$ an additive category with an $H$-action,
and let $\Gamma$ act on $\calc=\prod_{\Gamma/H}\calb$
as it does on the product in
Lemma~\ref{lem:on-homotopy-fixed-spectra}
~\ref{lem:on-h-f-spectra:coinduction}.
Then
\[
\IK^{-\infty}(\calc)^{\Gamma} \simeq
          \IK^{-\infty}(\calc)^{h_\calf \Gamma}.
\]
\end{lemma}

\begin{proof}
For $\calb=\cala(\pt)^\infty$ this is proven
in~\cite[Theorem 6.3]{Rosenthal-splitting-K-theory}
using Lemma~\ref{lem:on-homotopy-fixed-spectra}~
\ref{lem:on-h-f-spectra:coinduction}
and the fact that $\IK^{-\infty}$ commutes with
infinite products~\cite{Carlsson-K-theory-infinite-products}.
The same proof works in the general case.
\end{proof}

\begin{remark}
\label{rem:products-and-L-theory}
In order for the $L$-theory version of
Lemma~\ref{lem:products}
to be true, it must also be assumed that for sufficiently
large $i$, $K_{-i}(\calb[H])=0$.
This is needed because
the compatibility of $\IL^{-\infty}$ with infinite products
is only known if the $K$-theory vanishes in degree $-i$
for sufficiently large $i$ (see
\cite[p.756]{Carlsson-Pedersen-Controlled-algebra-Novikov}).
\end{remark}

Proposition~\ref{prop:vertical}, below,
is an important fact needed to prove the
Descent Principle (Theorem~\ref{thm:descent}).
Its proof is based on work of Carlsson and Pedersen
~\cite[Theorem 2.11]
{Carlsson-Pedersen-Controlled-algebra-Novikov}
who proved the result in the case
where $\Gamma$ is torsion-free and $X$ is a
finite $\Gamma$-CW complex
(meaning that $X$ has finitely many $\Gamma$-cells).
This was later generalized by the second
author~\cite{Rosenthal-splitting-K-theory}
to include groups with torsion.
The result proven here relaxes the finiteness condition on $X$,
requiring only that $X$ be a finite
dimensional $\Gamma$-CW complex.

\begin{proposition}\label{prop:vertical}
Let $\Gamma$ be a discrete group,
$\calf$ the family of finite subgroups
of $\Gamma$, $X$ a finite dimensional $\Gamma$-CW complex
with finite isotropy, and $p_{\Gamma}$ the quotient map
$X\to X/\Gamma$.
Then, for every $\Gamma$-invariant metric on $X$,
\[
{{\IK}}^{-\infty}
   \big( \cala_{p_{\Gamma}}(X)^\infty \big)^{\Gamma}
   \simeq
   {{\IK}}^{-\infty}\big(\cala_{p_{\Gamma}}(X)^\infty\big)^
                                           {h_\calf \Gamma}.
\]
\end{proposition}

\begin{proof}
Proceed by induction on the skeleta of $X$.
Let $X_0$ denote the 0-skeleton of $X$.
For some indexing set $J$, $X_0=\coprod_{j\in J}\Gamma/ H_j$,
where $H_j\in\calf$ for every $j\in J$.
Since we are taking germs away from zero and objects have
$\Gamma$-compact support,
\[
\cala_{p_{\Gamma}}(X_0)^\infty
    \cong \bigoplus_{j\in J}
    \Big(
      \prod_{\Gamma/H_j}\cala(\pt)^\infty
    \Big),
\]
which is a $\Gamma$-equivariant equivalence of categories.
Let
$\calc_j=\prod_{\Gamma/H_j}\cala(\pt)^\infty$.
Since $\IK^{-\infty}(\calc_j)$ is
$\Gamma$-invariant for every $j\in J$,
\[
\IK^{-\infty} \Big(\bigoplus_{j\in J} \calc_j\Big)^{\Gamma}
  \simeq \Big(\bigvee_{j\in J}
     \IK^{-\infty}(\calc_j)\Big)^{\Gamma}
  \simeq \bigvee_{j\in J} \IK^{-\infty}(\calc_j)^{\Gamma}
\]
and
\[
\IK^{-\infty} \Big(\bigoplus_{j\in J} \calc_j\Big)
                                    ^{h_\calf \Gamma}
   \simeq \Big(\bigvee_{j\in J}
       \IK^{-\infty}(\calc_j)\Big)^{h_\calf \Gamma}
   \simeq \bigvee_{j\in J}
       \IK^{-\infty}(\calc_j)^{h_\calf \Gamma},
\]
where the last weak equivalence makes use of
Lemma~\ref{lem:on-homotopy-fixed-spectra}
~\ref{lem:on-h-f-spectra:equivalence}.
By Lemma~\ref{lem:products},
\[
\IK^{-\infty}(\calc_j)^{\Gamma} \to
       \IK^{-\infty}(\calc_j)^{h_\calf \Gamma}
\]
is a weak homotopy equivalence.
Therefore,
\[
\bigvee_{j\in J} \IK^{-\infty}(\calc_j)^{\Gamma}
  \simeq \bigvee_{j\in J} \IK^{-\infty}(\calc_j)
  ^{h_\calf \Gamma},
\]
which completes the base case of the induction.

Now assume that the proposition holds for the
$(n-1)$-skeleton $X_{n-1}$.
Let $A \to B \to C$ denote the sequence
\[
{\IK}^{-\infty}\big(\cala_{p_\Gamma}(X_{n-1})^\infty\big)
  \to {\IK}^{-\infty}\big(\cala_{p_\Gamma}(X_n)^\infty\big)
  \to {\IK}^{-\infty}\big(\cala_{p_\Gamma}(X_n,X_{n-1})\big).
\]
Consider the following commutative diagram:
\[ \xymatrix{
    A^{\Gamma} \ar[d]_a \ar[r] &
    B^{\Gamma} \ar[d]_b \ar[r] &
    C^{\Gamma} \ar[d]_c
    \\
    A^{h_\calf \Gamma} \ar[r] &
    B^{h_\calf \Gamma} \ar[r] &
    C^{h_\calf \Gamma}. }
\]
Notice that each row in the diagram is a fibration of spectra.
The second row is a fibration since taking homotopy
fixed sets and taking homotopy fibers are both homotopy
limits and therefore commute~\cite{Bousfield-Kan}.

We must show that $b$ is a weak homotopy equivalence.
By the induction hypothesis, $a$ is a weak homotopy equivalence.
Therefore, by the Five Lemma, it suffices to prove that
$c$ is a weak homotopy equivalence.

Since we are taking germs away from $X^{n-1}$
(and thus away from zero)
and objects have
$\Gamma$-compact support,
\[
\cala_{p_\Gamma}(X_n,X_{n-1}) \cong
   \bigoplus_{i\in I}
          \Big(\prod_{\Gamma/H_i}\cala(D^n,S^{n-1})\Big).
\]
The proof is now completed by arguing as in the beginning
of the induction with
$\calc_i=\prod_{\Gamma/H_i}\cala(D^n,S^{n-1})$.
\end{proof}

\begin{theorem}[The Descent Principle]\label{thm:descent}
Let $\Gamma$ be a discrete group,
$\calf$ the family of finite subgroups of $\Gamma$,
$X$ a finite dimensional $\Gamma$-CW complex, and
$p_{\Gamma}$ the quotient map $X\to X/\Gamma$.
Assume that $X$ admits a $\Gamma$-invariant metric
such that $K_n\big(\cala^G_{p_\Gamma}(X)\big)=0$ for every
integer $n$, and every $G \in \calf$.

Then the map
$H_*^{\Gamma} (X;\bfK_\cala) \to K_*(\cala[\Gamma])$
is a split injection.
\end{theorem}

\begin{proof}
Consider the following commutative diagram of
fibration sequences:
\[
\xymatrix{
    \IK^{-\infty}\big(\cala_{p_\Gamma}(X)_0\big)^{\Gamma}
       \ar[d]_a \ar[r] &
    \IK^{-\infty}\big(\cala_{p_\Gamma}(X)\big)^{\Gamma}
       \ar[d]_b \ar[r] &
    \IK^{-\infty}\big(\cala_{p_\Gamma}(X)^\infty\big)^{\Gamma}
    \ar[d]_c
    \\
    \IK^{-\infty}\big(\cala_{p_\Gamma}(X)_0\big)^{h_\calf \Gamma}
      \ar[r] &
    \IK^{-\infty}\big(\cala_{p_\Gamma}(X)\big)^{h_\calf \Gamma}
      \ar[r] &
    \IK^{-\infty}\big(\cala_{p_\Gamma}(X)^\infty\big)^
                                         {h_\calf \Gamma}.
}
\]
By Proposition~\ref{prop:vertical},
$c$ is a weak homotopy equivalence.
Since each row in the diagram is a fibration,
it suffices to show that
$\IK^{-\infty}\big(\cala_{p_\Gamma}(X)\big)^{h_\calf \Gamma}$
is weakly contractible.
But this follows from
Lemma~\ref{lem:on-homotopy-fixed-spectra}~
\ref{lem:on-h-f-spectra:equivalence}
and the assumption that
$K_n\big(\cala^G_{p_\Gamma}(X)\big)=0$ for every
integer $n$.
\end{proof}


\section{Proof of the Main Theorem}
\label{sec:proof}

\begin{theorem}
\label{thm:split-injectivity-with-coeficients}
Let $\Gamma$ be a discrete group and let $\cala$ be a small
additive category.
Assume that there is a finite dimensional $\Gamma$-CW model
for the universal space for proper $\Gamma$-actions,
$\Ebar$, and assume that there is a
$\Gamma$-invariant metric on $\Ebar$ such that
$\Ebar$ is a complete proper path metric space
that is uniformly $\Fin$-contractible
and has finite asymptotic dimension.
Then the assembly map,
$H_*^{\Gamma} (\Ebar;
  \mathbb{K}^{-\infty}_\cala) \to K_*(\cala[\Gamma])$,
in algebraic $K$-theory, is a split injection.
\end{theorem}
Some preparations must be made before we can
give the proof of
Theorem~\ref{thm:split-injectivity-with-coeficients}.

Let $G$ be a finite subgroup of $\Gamma$.
Let
\[
G = H_0,\,H_1,\,\dots\,,\,H_m = \{e\}
\]
contain exactly one subgroup from each conjugacy class
of subgroups of $G$ and
let the $H_i$ be ordered by cardinality.
That is, $|H_i| \geq |H_{i+1}|$.

For each $k$, $0 \leq k \leq m$, define
$\cals_k = \big\{ H_i^g \; \big| \; 0 \leq i \leq k,\;
                                              g \in G \big\}$
and $Z_k = \Ebar^{\cals_k}$.
Clearly, $\cals_k$ is invariant under conjugation by $G$.
Therefore, $Z_k$ is $G$-invariant for every $k$,
$0 \leq k \leq m$.

\begin{notation}\label{notation:maps}
For each $k$,
let $p \colon Z_k \to \Ebar / \Gamma$ be the restriction
of the quotient map $p_{\Gamma}:\Ebar \to \Ebar / \Gamma$,
and let $p_G \colon Z_k / G \to \Ebar / \Gamma$ denote the
restriction of the canonical projection
$\Ebar / G \to \Ebar / \Gamma$.
\end{notation}

\begin{lemma}
\label{lem:K-theory-of-Z-k-mod-G-vanishes}
For every $k$, $0 \leq k \leq m$, and every subgroup $H$ of $G$,
the $K$-theory of $\cala[H]_{p_G}(Z_k / G)$ vanishes.
\end{lemma}

\begin{proof}
By Lemmas~\ref{lem:uniform-contractible-quotients}
and~\ref{lem:asymptotic-dimension-and-finite-quotients},
$Z_k / G$
is uniformly contractible and has finite asymptotic dimension.
Let $K$ be a subset of $Z_k / G$ whose image, $\bar K$,
in $\Ebar / G$ is compact.
If $\alpha > 0$, the image of $K^\alpha$ is $\bar K^\alpha$,
which
is compact since $\Ebar / G$ is a proper metric space.
Therefore, we can apply Proposition~\ref{prop:vanish}
to conclude that the $K$-theory of $\cala[H]_{p_G}(Z_k / G)$
vanishes.
\end{proof}

For the following fact,
compare~\cite[Lemma 7.4]{Rosenthal-splitting-K-theory}.

\begin{lemma}\label{lem:quotient}
For each $k$, $1\leq k\leq m$,
\[ \cala^G_p(Z_k,Z_{k-1}) \cong
            \cala[H_k]_{p_G}(Z_k/G,Z_{k-1}/G). \]
\end{lemma}

\begin{proof}
Since we are taking germs away from
$Z_{k-1}$ and $Z_{k-1} / G$,
every morphism has a representative that is
zero on $Z_{k-1}\times [0,1)$ and on
$Z_{k-1} / G\times [0,1)$, respectively.
Therefore, it is irrelevant what the objects over
$Z_{k-1}\times [0,1)$ and $Z_{k-1} / G\times [0,1)$ are.
By construction, the stabilizer subgroup of any point
not in $Z_{k-1}\times [0,1)$ is a conjugate of $H_k$.
Since $\cala^G_p(Z_k,Z_{k-1})$ is a fixed category,
the parts of an object over points in the same orbit
must be isomorphic modules.
Thus, the object $M$ in $\cala^G_p(Z_k,Z_{k-1})$
is sent to $M'$, where
$M'_{(y,t)}$ (with $(y,t)\notin Z_{k-1} / G\times [0,1)$)
is the $\cala[H_k]$-module sitting over the point
in $\big\{(x,t)\,\big|\,x\in p^{-1}(y)\big\}$
whose stabilizer subgroup is $H_k$.
The inverse of this is to take the
$\cala[H_k]$-module over $(y,t)$ and use the $G$-action to
spread it around the orbit
$\big\{(x,t)\,\big|\,x\in p^{-1}(y)\big\}$.
This explains how the objects in
$\cala^G_p(Z_k,Z_{k-1})$
and $\cala[H_k]_{p_G}(Z_k/G,Z_{k-1}/G)$
are identified.
To verify~\eqref{define-calx(Y,p)-compact},
notice that if $K_0$ is a subset of $Z_k / \Gamma$,
then the image of $p_{\Gamma}^{-1}(K_0)$
under the quotient map $Z_k\to Z_k / G$ is $p^{-1}(K_0)$.
Since we are taking germs, the components of a morphism need
to become small.
Therefore, non-zero components of a morphism have
the same isotropy, namely a conjugate of $H_k$.
Furthermore, the equivariance of morphisms in
$\cala^G_p(Z_k,Z_{k-1})$
implies that there is only one choice when lifting a
morphism from
$\cala[H_k]_{p_G}(Z_k/G,Z_{k-1}/G)$.
\end{proof}

\begin{proof}
[Proof of Theorem~\ref{thm:split-injectivity-with-coeficients}]
By the Descent Principle,
it suffices to show that the spectrum
$\IK^{-\infty}\big(\cala^G_{p_\Gamma}(\Ebar)\big)$
is weakly contractible for every finite subgroup $G$ of
$\Gamma$.

Let $G$ be a finite subgroup of $\Gamma$
and proceed by induction on the filtration
\[
\Ebar^G=Z_0 \subseteq Z_1 \subseteq
  \dddot{} \subseteq Z_{m-1} \subseteq Z_{m}=\Ebar
\]
defined above. Since $G$ acts trivially on $\Ebar^G$,
$\cala^G_p(\Ebar^G)$
is equivalent to
$\cala[G]_p(\Ebar^G)$.
By Lemma~\ref{lem:K-theory-of-Z-k-mod-G-vanishes},
$\IK^{-\infty}\big(\cala[G]_p(\Ebar^G)\big)$
is weakly contractible.
This completes the base case of the induction.

Assume now that $\IK^{-\infty}\big(\cala^G_p(Z_{k-1})\big)$
is weakly contractible.
We must show that $\IK^{-\infty}\big(\cala^G_p(Z_k)\big)$
is weakly contractible.
Consider the following Karoubi filtration
\[
\cala^G_p(Z_{k-1}) \to
 \cala^G_p(Z_k) \to
 \cala^G_p(Z_k,Z_{k-1})
\]
(see \eqref{pair-sequence-for-cala}), which yields a
fibration of spectra after applying
$\IK^{-\infty}$.
By using the induction hypothesis, we need only show that
$\IK^{-\infty}\big(\cala^G_p(Z_k,Z_{k-1})\big)$
is weakly contractible. By Lemma~\ref{lem:quotient},
$\cala^G_p(Z_k,Z_{k-1})$
is equivalent to $\cala[H_k]_{p_G}(Z_k/G,Z_{k-1}/G)$,
which fits into the Karoubi filtration:
\[
\cala[H_k]_{p_G}(Z_{k-1}/G) \to
     \cala[H_k]_{p_G}(Z_k/G) \to
     \cala[H_k]_{p_G}(Z_k/G,Z_{k-1}/G).
\]
Both $\IK^{-\infty}\big(\cala[H_k]_{p_G}(Z_{k-1}/G)\big)$
and $\IK^{-\infty}\big(\cala[H_k]_{p_G}(Z_k/G)\big)$
are weakly contractible by
Lemma~\ref{lem:K-theory-of-Z-k-mod-G-vanishes}.
Therefore,
$\IK^{-\infty}\big(\cala[H_k]_{p_G}(Z_k/G,Z_{k-1}/G)\big)$
is also weakly contractible.
\end{proof}


\section{$L$-theory}\label{sec:L-theory}

If $\cala$ is an additive category with involution
and an action of $\Gamma$, then there is an
assembly map
\[
H_*^\Gamma(\Ebar; \bfL_\cala) \to
H_*^\Gamma(\pt;\bfL_\cala)
= L_*^{\langle -\infty \rangle}(\cala[\Gamma])
\]
(see \cite[Section~5]{BR-coefficients}).
Here $\bfL_\cala$ is an $\Or (\Gamma)$-spectrum whose
value on $\Gamma / H$ is weakly equivalent
to $\IL^{-\infty} (\cala[H])$,
where $\IL^{-\infty}(\cala[H])$ is the spectrum
whose homotopy groups are
the {\em ultimate lower quadratic} $L$-groups
$L^{\langle -\infty \rangle}_*(\cala[H])$
(see \cite[Chapter 17]{Ranicki-lower}).
If $\cala$ is the category of finitely generated free
$R$-modules for a ring, $R$, with involution,
then $L^{\langle -\infty \rangle}_*(\cala[H]) =
      L^{\langle -\infty \rangle}_*(R[H])$,
and the above assembly map is
\[
H_*^\Gamma(\Ebar; \bfL_R) \to
H_*^\Gamma(\pt;\bfL_R) =
         L^{\langle -\infty \rangle}_*(R[\Gamma]).
\]
The following is the $L$-theory version of
Theorem~\ref{thm:split-injectivity-with-coeficients}.

\begin{theorem}
\label{thm:split-injectivity-with-coeficients-L-theory}
Let $\Gamma$ be discrete group and let $\cala$ be a small
additive category with involution.
Assume that there is a finite dimensional $\Gamma$-$CW$-model
for the universal space for proper $\Gamma$-actions,
$\Ebar$, and assume that there is a
$\Gamma$-invariant metric on $\Ebar$ such that
$\Ebar$ is uniformly $\Fin$-contractible,
is a complete proper path metric space
and has finite asymptotic dimension.
Assume that for each finite subgroup $G$
there is an $i_0 \in \IN$ such that for
$i \geq i_0$, $K_{-i}(\cala[G]) = 0$, where
the involution is forgotten and $\cala$ is considered
only as an additive
category.
Then the assembly map,
$H_*^{\Gamma} (\Ebar;
  \bfL_\cala) \to L_*(\cala[\Gamma])$,
in $L$-theory, is a split injection.
\end{theorem}

\begin{proof}
Everything we did for $K$-theory also works for
$L$-theory with the exception of Lemma~\ref{lem:products}.
As pointed out in Remark~\ref{rem:products-and-L-theory},
this lemma will carry over to $L$-theory if the
additional assumption about the vanishing
of $K_{-i}(\cala[G])$ for large $i$ is made.
The rest of the argument proceeds with out further changes.
For the required properties of $L$-theory, see
\cite[Section~4]{Carlsson-Pedersen-Controlled-algebra-Novikov}.
\end{proof}



\end{bibunit}


\begin{thebibliography}{ACFP94}


\bibitem[ACFP94]{ACFP-continuous-control}
D.~R. Anderson, F.~X. Connolly, S.~C. Ferry, and E.~K. Pedersen.
\newblock Algebraic {$K$}-theory with continuous control at infinity.
\newblock {\em J. Pure Appl. Algebra}, 94(1):25--47, 1994.

\bibitem[Bar03a]{Bartels-on-the-domain}
A.~Bartels.
\newblock On the domain of the assembly map in algebraic {$K$}-theory.
\newblock {\em Algebr. Geom. Topol.}, 3:1037--1050 (electronic), 2003.

\bibitem[Bar03b]{Bartels-finite-asymp-dim}
A.~Bartels.
\newblock Squeezing and higher algebraic {$K$}-theory.
\newblock {\em $K$-Theory}, 28(1):19--37, 2003.

\bibitem[BHM93]{Boekstedt-Hsiang-Madsen-cyclotomic-trace}
M.~B{\"o}kstedt, W.~C. Hsiang, and I.~Madsen.
\newblock The cyclotomic trace and algebraic ${K}$-theory of spaces.
\newblock {\em Invent. Math.}, 111(3):465--539, 1993.

\bibitem[BK72]{Bousfield-Kan}
A.K. Bousfield and D.M. Kan.
\newblock {\em Homotopy Limits, Completions and Localizations}.
\newblock Springer-Verlag, New York, 1972.
\newblock Lecture Notes in Mathematics, no. 304.

\bibitem[BR05]{BR-coefficients}
A.~Bartels and H.~Reich.
\newblock Coefficients for the {F}arrell-{J}ones conjecture.
\newblock Preprintreihe SFB 478 --- Geometrische Strukturen in der Mathematik,
  Heft 402, M\"unster, arXiv:math.KT/0510602, 2005.

\bibitem[Car95]{Carlsson-K-theory-infinite-products}
G.~Carlsson.
\newblock On the algebraic {$K$}-theory of infinite product categories.
\newblock {\em $K$-Theory}, 9(4):305--322, 1995.

\bibitem[CG04a]{Carlsson-Goldfarb-finite-asymptotic}
G.~Carlsson and B.~Goldfarb.
\newblock The integral {$K$}-theoretic {N}ovikov conjecture for groups with
  finite asymptotic dimension.
\newblock {\em Invent. Math.}, 157(2):405--418, 2004.

\bibitem[CG04b]{Carlsson-Goldfarb-coherence}
G.~Carlsson and B.~Goldfarb.
\newblock On homological coherence of discrete groups.
\newblock {\em J. Algebra}, 276(2):502--514, 2004.

\bibitem[CP95]{Carlsson-Pedersen-Controlled-algebra-Novikov}
G.~Carlsson and E.~K. Pedersen.
\newblock Controlled algebra and the {N}ovikov conjectures for {$K$}- and
  {$L$}-theory.
\newblock {\em Topology}, 34(3):731--758, 1995.

\bibitem[DL98]{Davis-Lueck-assembly}
J.~F. Davis and W.~L{\"u}ck.
\newblock Spaces over a category and assembly maps in isomorphism conjectures
  in {$K$}- and {$L$}-theory.
\newblock {\em $K$-Theory}, 15(3):201--252, 1998.

\bibitem[FJ93]{Farrell-Jones-Isomorphism-Conjectures}
F.~T. Farrell and L.~E. Jones.
\newblock Isomorphism conjectures in algebraic {$K$}-theory.
\newblock {\em J. Amer. Math. Soc.}, 6(2):249--297, 1993.

\bibitem[FW91]{Ferry-Weinberger-Cuvature-tangentiality+controlled}
S.~C. Ferry and S.~Weinberger.
\newblock Curvature, tangentiality, and controlled topology.
\newblock {\em Invent. Math.}, 105(2):401--414, 1991.

\bibitem[Gro93]{Gromov-asymptotic-invariants}
M.~Gromov.
\newblock Asymptotic invariants of infinite groups.
\newblock In {\em Geometric group theory, Vol.\ 2 (Sussex, 1991)}, pages
  1--295. Cambridge Univ. Press, Cambridge, 1993.

\bibitem[Hel62]{Helgason-Book-1962}
S.~Helgason.
\newblock {\em Differential geometry and symmetric spaces}.
\newblock Pure and Applied Mathematics, Vol. XII. Academic Press, New York,
  1962.

\bibitem[Hig00]{Higson-Bivariant-K-theory}
N.~Higson.
\newblock Bivariant {$K$}-theory and the novikov conjecture.
\newblock {\em Geom. Funct. Anal.}, 10(3):563--581, 2000.

\bibitem[Hoc65]{Hochschild-Lie-groups-book}
G.~Hochschild.
\newblock {\em The structure of {L}ie groups}.
\newblock Holden-Day Inc., San Francisco, 1965.

\bibitem[HP04]{Hambleton-Pedersen-identify}
I.~Hambleton and E.~K. Pedersen.
\newblock Identifying assembly maps in {$K$}- and {$L$}-theory.
\newblock {\em Math. Ann.}, 328(1-2):27--57, 2004.

\bibitem[HR95]{Higson-Roe-On-coarse-Baum-Connes}
N.~Higson and J.~Roe.
\newblock On the coarse {B}aum-{C}onnes conjecture.
\newblock In {\em Novikov conjectures, index theorems and rigidity, Vol.\ 2
  (Oberwolfach, 1993)}, volume 227 of {\em London Math. Soc. Lecture Note
  Ser.}, pages 227--254. Cambridge Univ. Press, Cambridge, 1995.

\bibitem[HR00]{Higson-Roe-amenable-actions+Novikov}
N.~Higson and J.~Roe.
\newblock Amenable group actions and the {N}ovikov conjecture.
\newblock {\em J. Reine Angew. Math.}, 519:143--153, 2000.

\bibitem[Ji04]{Ji-asymptotic-Novikov-arithmetic}
L.~Ji.
\newblock Asymptotic dimension and the integral {$K$}-theoretic {N}ovikov
  conjecture for arithmetic groups.
\newblock {\em J. Differential Geom.}, 68(3):535--544, 2004.

\bibitem[Kas88]{Kasparov-equivariant-KK-and-Novikov}
G.~G. Kasparov.
\newblock Equivariant {$KK$}-theory and the {N}ovikov conjecture.
\newblock {\em Invent. Math.}, 91(1):147--201, 1988.

\bibitem[LR05]{Lueck-Reich-survey}
W.~L{\"u}ck and H.~Reich.
\newblock The {B}aum-{C}onnes and the {F}arrell-{J}ones conjectures in {$K$}-
  and {$L$}-theory.
\newblock In {\em Handbook of $K$-theory. Vol. 1, 2}, pages 703--842. Springer,
  Berlin, 2005.

\bibitem[LRV03]{Lueck-Reich-Varisco-commuting}
W.~L{{\"u}}ck, H.~Reich, and M.~Varisco.
\newblock Commuting homotopy limits and smash products.
\newblock {\em $K$-Theory}, 30(2):137--165, 2003.
\newblock Special issue in honor of Hyman Bass on his seventieth birthday. Part
  II.

\bibitem[L{\"u}c02]{Lueck-Chern-character}
W.~L{\"u}ck.
\newblock Chern characters for proper equivariant homology theories and
  applications to {$K$}- and {$L$}-theory.
\newblock {\em J. Reine Angew. Math.}, 543:193--234, 2002.

\bibitem[L{\"u}c05]{Lueck-survey-classifying-spaces}
W.~L{\"u}ck.
\newblock Survey on classifying spaces for families of subgroups.
\newblock In {\em Infinite groups: geometric, combinatorial and dynamical
  aspects}, volume 248 of {\em Progr. Math.}, pages 269--322. Birkh\"auser,
  Basel, 2005.

\bibitem[PW85]{Pedersen-Weibel-delooping}
E.~K. Pedersen and C.~A. Weibel.
\newblock A nonconnective delooping of algebraic {$K$}-theory.
\newblock In {\em Algebraic and geometric topology (New Brunswick, N.J.,
  1983)}, volume 1126 of {\em Lecture Notes in Math.}, pages 166--181.
  Springer, Berlin, 1985.

\bibitem[Ran92]{Ranicki-lower}
A.~A. Ranicki.
\newblock {\em Lower {$K$}- and {$L$}-theory}, volume 178 of {\em London
  Mathematical Society Lecture Note Series}.
\newblock Cambridge University Press, Cambridge, 1992.

\bibitem[Roe91]{Roe-hyperbolic-spaces-exotic-cohomology}
J.~Roe.
\newblock Hyperbolic metric spaces and the exotic cohomology {N}ovikov
  conjecture.
\newblock {\em $K$-Theory}, 4(6):501--512, 1990/91.

\bibitem[Ros04]{Rosenthal-splitting-K-theory}
D.~Rosenthal.
\newblock Splitting with continuous control in algebraic {$K$}-theory.
\newblock {\em $K$-Theory}, 32(2):139--166, 2004.

\bibitem[Wei02]{Weiss-excision-and-restriction}
M.~Weiss.
\newblock Excision and restriction in controlled {$K$}-theory.
\newblock {\em Forum Math.}, 14(1):85--119, 2002.

\bibitem[Yu98]{Yu-finite-asymptotic-dimension-Novikov}
G.~Yu.
\newblock The {N}ovikov conjecture for groups with finite asymptotic dimension.
\newblock {\em Ann. of Math. (2)}, 147(2):325--355, 1998.

\end{thebibliography}

\begin{thebibliography}{1}

\bibitem{Bartels-Rosenthal-On-K-fad}
A.~Bartels and D.~Rosenthal.
\newblock On the {$K$}-theory of groups with finite asymptotic dimension.
\newblock {\em J. Reine Angew. Math.}, 612:35--57, 2007.

\bibitem{Kasprowski-FDC}
D.~Kasprowski.
\newblock On the {$K$}-theory of groups with finite decomposition complexity.
\newblock {\em Proc. Lond. Math. Soc. (3)}, 110(3):565--592, 2015.

\bibitem{Kasprowski-discr-Lie}
D.~Kasprowski.
\newblock On the {$K$}-theory of subgroups of virtually connected {L}ie groups.
\newblock {\em Algebr. Geom. Topol.}, 15(6):3467--3483, 2015.

\bibitem{Rosenthal-splitting}
D.~Rosenthal.
\newblock Splitting with continuous control in algebraic {$K$}-theory.
\newblock {\em $K$-Theory}, 32(2):139--166, 2004.

\end{thebibliography}

\parindent0em

\begin{bibunit}

\newpage

\thispagestyle{empty}

\begin{center}
{\bf \large Erratum for ``On the $K$-theory of groups with
                   finite asymptotic dimension"}\\
                   \vspace{3ex}
     Arthur Bartels and David Rosenthal              
\end{center}

\vspace{4ex}

Daniel Kasprowski pointed out an error in our treatment of the descent principle in~\cite{Bartels-Rosenthal-On-K-fad}.
  The problem is that the map
  \[
   \Big(\bigvee_{j\in J}
       \IK^{-\infty}(\calc_j)\Big)^{h_\calf \Gamma}
   \to \bigvee_{j\in J}
       \IK^{-\infty}(\calc_j)^{h_\calf \Gamma},
  \]
  in the proof of Proposition~7.5 is only an equivalence if there exists a cocompact model for the classifying space for proper actions ${\underbar{\rm{E}}\Gamma}$.
  This means that we did not prove a more general descent principle than~\cite{Rosenthal-splitting} and that Theorems~A and~B are only proved under the assumption that there exists a cocompact model for ${\underbar{\rm{E}}\Gamma}$. 
  This similarly affects the two corollaries formulated in the introduction; i.e., our proof of split  injectivity for the assembly maps in $K$- and $L$-theory relative to the family of finite subgroups is only valid for discrete subgroups of virtually connected Lie groups that admit a cocompact model for ${\underbar{\rm{E}}\Gamma}$. 

  Daniel Kasprowski developed an alternative descent principle using the concept of \emph{bounded homotopy fixed points} and proved versions of Theorems~A and~B under the assumption that there exists a finite dimensional model for ${\underbar{\rm{E}}\Gamma}$ and a bound on the order of the finite subgroups of $\Gamma$~\cite{Kasprowski-FDC}.
  Concerning discrete subgroups of Lie groups, he proved split injectivity of the assembly maps in question for all finitely generated discrete subgroups of virtually connected Lie groups~\cite{Kasprowski-discr-Lie}.

\end{bibunit}

%

%


\end{document}